\documentclass[11pt]{amsart}
\usepackage{amsmath,amssymb, graphics, amscd,latexsym,here}
\makeatletter

\newtheorem{Theorem}{Theorem}
\newtheorem{Lemma}[Theorem]{Lemma}
\newtheorem{Corollary}[Theorem]{Corollary}
\newtheorem{Proposition}[Theorem]{Proposition}

\newtheorem{Remark}[Theorem]{Remark}

\newcommand{\eps}{\varepsilon}

\newcommand\la{\lambda}

\newcommand\al{\alpha}
\newcommand\La{\Lambda}

\newcommand\be{\beta}
\newcommand\Si{\Sigma}
\newcommand\ga{\gamma}

\newcommand\de{\delta}
\newcommand\De{\Delta}

\newcommand\cM{\mathcal  M}

\newcommand\cP{\mathcal  P}

\newcommand\cS{\mathcal  S}

\newcommand\BC{ {\Bbb C}}

\newcommand\BZ{{\Bbb  Z}}

\newcommand\BP{ {\Bbb  P}}

\newcommand\bfN{\mbox {\bf  N}}

\newcommand\bfw{\mbox {\bf  w}}

\newcommand\nl{\newline}

\newcommand\wtl{\widetilde}

\newcommand\tAS{\text{t-AS}\/}

\newcommand\Ker{\text{Ker}\/}

\newcommand\id{\text{id}}

\newcommand\degree{\text{degree}\/}

\newcommand\inv{^{-1}}

\def\mapright#1{\smash{\mathop{\longrightarrow}\limits^{{#1}}}}

\def\mapdown#1{\Big\downarrow\rlap{$\vcenter{\hbox{$#1$}}$}}

\def\inv{^{-1}}

\begin{document}
\title[Tangential Alexander polynomials 
]
{Tangential Alexander polynomials 
and non-reduced degeneration
}

\author
[M. Oka ]
{Mutsuo Oka }
\address{\vtop{
\hbox{Department of Mathematics}
\hbox{Tokyo  University of Science}
\hbox{26 Wakamiya-cho, Shinjuku-ku}
\hbox{Tokyo 162-8601}
\hbox{\rm{E-mail}: {\rm oka@rs.kagu.tus.ac.jp}}
}}

\begin{abstract}
We introduce a notion of tangential Alexander polynomials
for plane curves
 and study the relation
with $\theta$-Alexander polynomial.  As an application,
  we use these polynomials to study
a  non-reduced degeneration  $C_t,\,\to \,D_0+jL$.
 We show that
there exists a certain surjectivity of the fundamental groups
and  divisibility among their Alexander polynomials.
\end{abstract}
\maketitle

\maketitle

\section{Introduction}


Let $C$ be a plane curve.
We are interested in the geometry of plane curves. Choose a line
$L\subset \BP^2$  and put  $\BC_L^2:=\BP^2-L$.
As geometrical invariants, we consider

(a) Fundamental groups: $\pi_1(\BP^2-C)$
and  $\pi_1(\BC_L^2-C)$

(b) Alexander polynomial $\De_C(t;L)$.

 \noindent
 Zariski studied $\pi_1(\BP^2-C)$ systematically
\cite{Za1}
and further developments have been made by many authors. 
%
To compute the Alexander polynomial, we need to choose a line at
infinity $L$.
However for a generic $L$, the Alexander polynomial has too much restrictions and
 we have often the trivial case $\De_C(t;L)=(t-1)^{r-1}$
where $r$ is the number of the irreducible components.
In our previous paper \cite{OkaSurvey}, we have introduced 
the notion of $\theta$-Alexander polynomials. This gives more
 informations 
for certain reducible curves but it  does not give any further 
information for irreducible curves.

 The purpose of this paper is to
 introduce the notion of the tangential Alexander polynomials.
Namely we consider all tangent lines $T_P C$ for the line at infinity.
It turns out that tangential Alexander polynomials
are  related to $\theta$-Alexander polynomials. We
apply these polynomials to study  non-reduced degenerations.
Let $C_t,\,t\in \De$ be a analytic  family of reduced curves for $t\ne
 0$
such that $C_t\,\mapright{t\to \infty}\,C_0=D_0+jL$ where
 $L$ is a line. The case $j\ge 2$ is a typical
 {\em non-reduced degeneration}. In this situation we 
study the geometry of $D_0$ using that of $C_t$. One of our results is 
the surjectivity assertion of the natural homomorphism:
$$ \phi:\,\pi_1(\BC_L^2-D_0) \to \pi_1(\BC_L^2-C_t)$$
Here  the point is that  $L$ is the line component of the limit curve $C_0$
( Theorem \ref{surjectivity}, \S 5).
This paper consists of the following sections.
\newline
\S 2 Fundamental groups
\newline
\S 3  Alexander polynomial
\newline
\S 4  {Dual stratification and tangential fundamental groups}
\newline
\S  5  Degeneration into non-reduced curves with a multiple line

\section{Fundamental groups}
Let $L$ be a fixed line and put $\BC_L^2:=\BP^2-L$.
We say {\em $L$ is generic with respect to $C$} if $L$ and $C$ intersects
transversely. The topology of $\BC_L^2-C$
does not depend on $L$ if $L$ is generic  and we call it 
{\em  the generic affine
complement}
and we often write as $\BC^2$ instead of $\BC_L^2$.
The following Lemma describes the relation of two fundamental groups.
\begin{Lemma}\label{Affine-Projective}{\rm (\cite{OkaOrdinary})} 
 Let $\omega$ be
 a lasso for $L$ and $N(\omega)$ be the subgroup normally generated by 
$\omega$. 

\noindent
{\rm (1)}  The  following  sequence is  exact.
$$1\to \bfN(\omega)\to \pi_1(\BC_L^2- C,b_0)
\to \pi_1(\BP^2-C,b_0)\to 1$$

\noindent
{\rm (2)} Assume that   $L$ is  generic. Then

{\rm (i)} $\omega$ is in the center of $\pi_1(\BC_L^2-C)$ and 
$\bfN(\omega)\cong \BZ$.

{\rm{(ii)}} We have the equality 
$D(\pi_1(\BC^2-C))=D(\pi_1(\BP^2-C))$
among their commutator groups.
 Thus $\pi_1(\BP^2-C)$ is
 abelian  if and only if $\pi_1(\BC^2-C)$ is abelian.
\end{Lemma}
 For non-generic line $L$, $\pi_1(\BC_L^2-C)$ may be
 non-abelian even if $\pi_1(\BP^2-C)$ is abelian.
For example, let $C=\{Y^2Z-X^3=0\}$ and take $L=\{Z=0\}$.
Then $\pi_1(\BC_L^2-C)\cong B_3$
where $B_3$ is the braid group of three strings and we recall that 
$B_3\cong \langle a,b\,|\, aba=bab\rangle$ (\cite{Artin1}).

%
\subsection {\bf First homology group $H_1(\BP^2-C)$.} 
Assume that
$C$ is a projective curve with $r$ irreducible components $C_1,\dots, C_r$
of degree $d_1,\dots, d_r$ respectively. 
By Lefschetz duality,
we have the following.
\begin{Proposition}  $H_1(\BP^2-C,\BZ)$
 is isomorphic to
$\BZ^{r-1}\times (\BZ/d_0\BZ)$ where  $d_0=\gcd (d_1,\dots,d_r)$. 
 In particular,  $H_1(\BC_L^2-C)\cong \BZ^r$.
\end{Proposition}
Take a lasso $g_i$ for each  component $C_i$ of $C$ for $i=1,\dots, r$. 
Then the corresponding homology classes 
$\{[g_i],\,i=1,\dots, r\}$  give free abelian generators
of $ H_1(\BC_L^2-C)$.
\subsection {\bf Degenerations and fundamental groups.}
Let $C$ be a reduced plane curve. The {\em total Milnor number} $\mu(C)$
is defined by the sum of the local Milnor numbers $\mu(C,P)$ 
at  the singular points $P$ of $ C$. Let $\De:=\{\zeta\in \BC\,|\,|\zeta|\le
1\}$
the unit disk.
We consider an analytic
 family of projective curves $C_t=\{F_t(X,Y,Z)=0\},\,t\in \De$
where
$F_t(X,Y,Z)$ are reduced  homogeneous polynomial of degree $d$ for any
$t$.
We call $\{C_t;t\in \De\}$ {\em a reduced degeneration}.
 We assume
that $C_t,~t\ne 0$ have the same configuration of
singularities so that they  are topologically equivalent but $C_0$  may obtain more singularities,
i.e., $\mu(C_t)\le\mu(C_0)$.
\begin{Theorem}\label{surjectivity-reduced}
For a given reduced degeneration $\{C_t;t\in \De\}$,
there is a canonical surjective homomorphism
 for  $t\ne 0$:
\[
 \varphi:\pi_1(\BP^2-C_0)\twoheadrightarrow \pi_1(\BP^2-C_t)
\]
In particular, if $\pi_1(\BP^2-C_0)$ is abelian, so is $ \pi_1(\BP^2-C_t)$.
\end{Theorem}
See for example, \cite{OkaSurvey}
and also Theorem \ref{surjectivity} of \S 5 for  another simple proof.

\subsection{Product formula.}
 Assume that $C_i$ is a curve of degree $d_i,\,i=1,2$
which 
are intersecting  transversely at $d_1d_2$ distinct
points.
We denote the transversality  as  $C_1\pitchfork C_2$.
Take a line $L$ such that  $L\cap C_1\cap C_2=\emptyset$.
Note that $L$ need not be generic for $C_1$ or $C_2$.
\begin{Theorem}\label{Oka-Sakamoto}{\rm
 (Oka-Sakamoto\,\cite{Oka-Sakamoto})}
Under the above assumption, we have
$$
\pi_1(\BC_L^2-C_1\cup C_2)\cong\pi_1(\BC_L^2-C_1)\times\pi_1(\BC_L^2-C_2)$$
 \end{Theorem}
\noindent
For further information about fundamental groups, we refer to 
\cite{Artin1,Libgober-survey,OkaCertain,OkaTwo,Shimada1}.

\subsection{Example} \subsubsection{Abelian cases}
A curve $C$ with small singularities has often commutative fundamental group
$\pi_1(\BP^2-C)$.
Some examples are here:

\noindent
-- $C$ is a smooth irreducible curve.

\noindent
-- Irreducible curves with only $A_1$-singularities (i.e., nodes) by
 \cite{Za1,Fulton2, Deligne2, Harris1,Nori1,Orevkov1},
or irreducible curve of degree $d$ with $a$ nodes and $b$
cusps (i.e., $A_2$) with $6b+2a<d^2$ (\cite{Nori1}).

\noindent
-- $\pi_1(\BP^2-C)$ (respectively $\pi_1(\BC_L^2-C)$)
is abelian for  any irreducible curve of degree $d$
if it has a flex of order $\ge d -3$ in $\BC_L^2$ (resp. of order $d-2$) (\cite{Za1}).

 Let $f:\BC^2\to \BC$ be a polynomial mapping. Recall
 that $\alpha$ is a {\em atypical value at infinity} if the topological 
triviality at infinity fails at $t=\alpha$ for the family
of curves $C_t:=\overline{f\inv (t)}$
(see \cite{Ha-Le}).  
\begin{Proposition}\label{Prop6.2}{\rm (\cite{OkaTwo})}
 Let $f:\BC^2\to \BC$ be a polynomial mapping and
 assume that $0$ is not an { atypical value at infinity}
 and $C=f\inv(0)$ is smooth
in $\BC^2$. Then 
$\pi_1(\BC^2-C)\cong \BZ$.
\end{Proposition}

\subsubsection{Non-abelian case}
Assume that $p,q$ are positive integers greater than 1 and consider
the curve
\[
 C_{p,q}: \quad f_p(X,Y,Z)^q+f_q(X,Y,Z)^p=0
\]
where $f_p,\, f_q$ are  polynomials of degree $p,\,q$ respectively.
$C_{p,q}$ is called {\em a curve of $(p,q)$-torus type}.
Assume that two curves $\{f_p=0\}$ and $ \{f_q=0\}$ intersect transversely
and there is no other singularities  of $C_{p.q}$. Then
$\pi_1(\BP^2-C_{p,q})\cong G(p,q,q)$ and 
$\pi_1(\BC^2-C)\cong G(p,q)$.
 In particular, 
if $p,\,q$ are coprime,
$\pi_1(\BP^2-C_{p,q})\cong \BZ_p*\BZ_q$.
For the definition of $G(p,q)$ and $G(p,q,r)$,
we refer to \cite{OkaCertain}.
  
\subsection{Class formula and flex formula}\label{flexformula}
For the study of curves of low degree, it is often important to know 
the existence  of flex points.
 Let $d=\degree\, (C)$, $\check d$ be
the 
{\em degree of the dual curve} $\check C$, let $\Si(C)$ be the
singular points of $C$ and let
$\alpha(C)$ be {\em the number of the 
flex points}. Then $\check d$ and $\alpha(C)$ are given by the formula:
\begin{eqnarray*}
&\check d=d(d-1)\,-\, \sum_{P\in \Si(C)}(\mu(C,P)+m(C,P)-1)\\
& \alpha(C)=3d(d-2)\,-\, \sum_{P\in \Si(C)} \gamma(C,P)
\end{eqnarray*}
where $m(C,p)$ is the multiplicity of $C$ at $P$ and
 $\gamma(C,P)$ is the flex defect of the singularity 
$(C,P)$ \cite{NambaBook,Okadual}. (In
\cite{Okadual}, we  have denoted $\gamma(C,P)$ as $\de(C,P)$.
To distinguish with $\de$-genus of the singularity, we change the notation.)
\section{Alexander polynomial}
\subsection{General definition.}
Let  $X$ be  a finite connected  CW-complex
and let
$ \phi: \pi_1(X)\to \BZ$
be a surjective homomorphism. We fix a generator $t$ of the infinite
cyclic group $\BZ$.
Let  $\La$ is the group ring of $\BZ$. Then $\La$
is isomorphic to the Laurent polynomial ring $\BC[t,t\inv]$ and $\La$ is
 a principal ideal domain. Consider an infinite cyclic
covering
 $p: \wtl X\to X$ 
such that  $p_{\#}(\pi_1(\wtl X))=\Ker\,\phi$.
Then  $H_1(\wtl X,\BC)$ has a structure of $\La$-module where $t$ acts as
the canonical covering transformation. Thus
by the structure theorem of modules over a principal ideal domain,
 we have an identification:
\[
 H_1(\wtl X,\BC)\cong \La/\la_1\oplus\cdots\oplus \La/\la_n
\]
as $\La$-modules.
We  normalize the denominators so that $\la_i$ is a polynomial in $t$
with $\la_i(0)\ne 0$ for each $i=1,\dots, n$.
{\em The Alexander polynomial associated to $\phi$} is defined (see
\cite{Libgober1}) by 
the product
$\De_\phi(t):=\prod_{i=1}^n \la_i(t).$
\subsection{Alexander polynomials of plane curves}In our situation, we consider
 a plane curve $C=C_1\cup \cdots\cup C_r$ 
where  $C_1,\dots, C_r$ are irreducible components of
degree $d_1,\dots, d_r$ respectively.
Take a line $L$ as the line at infinity and 
let $\phi_\theta$ be  the composition 
\[
 \phi_\theta:\, \pi_1(\BC_L^2-C)\mapright{\xi}H_1(\BC_L^2-C,\BZ)
\cong  \BZ^r\overset{\theta}\to \BZ
\]
where $\theta$ is a surjective homomorphism.
Recall that  $\theta$ is determined by giving  an integer $n_i$ to
 each component $C_i$
such that $\gcd(n_1,\dots, n_r)=1$. We call $n_i$ {\em the weight for
$C_i$}.
{\em The Alexander polynomial of $C$ with respect to $(L,\theta)$}
is defined by $\De_{\phi_\theta}(t)$ and we denote it as 
$\De_C(t;L,\theta)$. 

\noindent
(1)  (Generic case) Assume that 
  $L$ to be generic and $\theta=\theta_{sum}$ where
$\theta_{sum}$ is defined by the canonical summation
 $\theta_{sum} (a_1,\dots,a_r)= \sum_{i=1}^r a_i$
( weight 1 for each component.)
In this  case, we simply write as $\De_C(t)$ and call it
{\em the generic
Alexander polynomial of $C$}, 
as it does not depend on the choice of a generic $L$.

\noindent
(2) If $\theta$ is the canonical summation $\theta_{sum}$
 but $L$ is not generic, 
we denote it as $\De_C(t;L)$, omitting $\theta$. 
 In particular, when $L$ is the tangent
 line of a smooth point $P\in C$, we call 
$\De_C(t;L)$ {\em the tangential Alexander polynomial at $P$ } and we also
 use the notation $\De_C(t;P)$.

\noindent
(3) If $L$ is generic  but $\theta$ is not  $\theta_{sum}$,
 we called $\De_C(t;L,\theta)$ {\em  the
$\theta$-Alexander polynomial} and we denote it as
$\De_C(t;\theta)$. 
In \cite{OkaSurvey} we denoted it by  $\De_{C,\theta}(t)$,
 but for the consistency of the
 notation  with (2),
we use the notation $\De_C(t;\theta)$.

\vspace{.2cm}\noindent
Recall that  
$(t-1)^{r-1}\,|\, \De_C(t)$ (\cite{OkaSurvey}). Thus this is also the
case
for  $\De_C(t;L)$ with any line $L$,
as $\De_C(t)\,|\,\De_C(t;L)$. We say that
$\De_C(t:L)$ is {\em trivial} if $\De_C(t;L)=(t-1)^{r-1}$.
\subsection{Fox calculus.}
Suppose that $G$ is a group and
$\phi:G\to \BZ$ is a given surjective homomorphism.
Assume that  $G$ has  a  finite presentation
as
\begin{eqnarray*}
 G\cong \langle
x_1,\dots, x_n\,|\, R_1,\dots, R_m\rangle\end{eqnarray*}
This  corresponds to a surjective homomorphism $\psi:\,F(n)\to G$
so that $\Ker\, \psi$ is normally generated by the words $R_1,\dots, R_m$
 where $F(n)$
is a free group of rank $n$, generated by $x_1,\dots,x_n$.
Consider the group ring $\BC(F(n))$  of $F(n)$ with $\BC$-coefficients.
The {\em  Fox differentials}
$ \frac{\partial}{\partial x_j}:\BC(F(n))\to \BC(F(n))$ for $j=1,\dots,
n$,
are additive homomorphisms which 
are characterized by the following properties.
\begin{eqnarray*}
(1)\, \frac{\partial }{\partial x_j} x_i=\de_{i,j},\quad
(2)\, \text{(Leibniz rule)}\,\, \frac{\partial}{\partial x_j}(uv)=\frac{\partial u}{\partial x_j}+u
\frac{\partial v}{\partial x_j},\, u,v\in \BC(F(n))
\end{eqnarray*}
The composition $\phi\circ \psi:F(n)\to \BZ$
gives a ring homomorphism
$\gamma: \BC(F(n))\to \BC[t,t\inv]$.
The {\em  Alexander matrix} $A$ is an $m\times n$ matrix with coefficients in 
$\BC[t,t\inv]$ and its $(i,j)$-component is given by
$\gamma( \frac{\partial R_i}{\partial x_j})$.
Then the {\em  Alexander polynomial} $\De_\phi(t)$
is  defined  by the greatest common divisor of $(n-1)\times
(n-1)$-minors 
of $A$. In the case of $G=\pi_1(X)$ for some connected 
topological space $X$, this definition coincides with the previous one
 (Fox \cite{Fox}).

\subsubsection{Examples}We gives several examples.
\nl
1. Consider the trivial case: $G=\BZ^r$ and $\phi=\theta_{sum}$, the 
canonical  one. Then 

1-1. $G\cong \BZ\cong \langle x_1\rangle$. Then $ \De(t)=1$.

1-2. If $ G= \BZ^r\cong\langle x_1,\dots,x_r\,|\, R_{i,j}=x_ix_j{x_i}\inv {x_j}\inv ,\, 1\le i<j\le r\rangle$,
 we have $\De(t)=(t-1)^{r-1}$. This follows from the Fox derivation:
\[
 \frac{\partial}{\partial x_i}(x_ix_j{x_i}\inv
 {x_j}\inv)=1-x_ix_jx_i\inv,
\quad
 \frac{\partial}{\partial x_j}(x_ix_j{x_i}\inv
 {x_j}\inv)=x_i-x_ix_jx_i\inv x_j\inv.
\]

\noindent
2. Let $C=\{Y^2Z-X^3=0\}$ and $L_{gen}=\{Z=Y\},\,L=\{Z=0\}$.
Note that $(0,1,0)$ is a flex point of $C$ and $L$ is the flex tangent.
Then
\begin{eqnarray*}
&\pi_1(\BC_{L_{gen}}^2-C)\cong \BZ,\quad\De_C(t;L_{gen})=1\\
& \pi_1(\BC_L-C)=\langle x_1,x_2\,|\,x_1x_2x_1=x_2x_1x_2\rangle\cong B_3,
\quad\De_C(t;L)=t^2-t+1
\end{eqnarray*}

\noindent
3. Let us consider the curve $C= \{Y^2Z^3-X^5=0\}\subset \BC^2$ and
$ L=\{Z=0\},\,M=\{Y=0\}$.
 Then
$\pi_1(\BP^2-C)\cong \BZ/5\BZ$  and
$\pi_1(\BC_L^2-C)\cong G(2,5)$ and 
$ \pi_1(\BC_M^2-C)\cong G(3,5)$.
In this case, we get 
\[
\De_C(t)=1,\quad \De_C(t;L)
 =\frac{(t^{10}-1)(t-1)}{(t^5-1)(t^2-1)},\quad
\De_C(t;M)=\frac{(t^{15}-1)(t-1)}{(t^5-1)(t^3-1)}
\]

\subsection{Weakness of the generic Alexander polynomial $\De_C(t)$.}
The following  Lemma describes the relation between the  Alexander 
polynomial and local singularities.
\begin{Lemma}\label{localAlexander}{\rm (Libgober \cite{Libgober1})}
Let $P_1,\dots, P_k$ be the singular points of $C$ (including those at
 infinity) and let 
$\De_i(t)$ be the characteristic polynomial of the Milnor fibration 
of the germ $(C,P_i)$. Then the generic Alexander polynomial
satisfies the divisibility:
$\De_C(t;L)\,|\,\prod_{i=1}^k \De_i(t)$.
\end{Lemma}

\begin{Lemma}\label{AlexanderInfinity}{\rm (Libgober \cite{Libgober1})}
 Let $d$ be the degree of $C$. Then the Alexander polynomial
$\De_C(t;L_\infty)$ divides the Alexander polynomial at infinity
$\De_\infty(t)$.
If $L_\infty$ is generic, $\De_\infty(t)$
 is given by 
$(t^d-1)^{d-2}(t-1)$.
In particular,  the roots of the generic 
Alexander polynomial are $d$-th roots of unity.
\end{Lemma}
\begin{Corollary}{\rm (\cite{Libgober1}, See also \cite{Za1})} Assume that $C$ is an irreducible
 curve
of degree $d$ and assume 
that the singularities  are either nodes (i.e., $A_1$) or ordinary cusp 
singularities (i.e., $A_2$). If $d$
is not divisible by $6$, the generic Alexander polynomial $\De_C(t)$ is
trivial.
\end{Corollary} 
This implies that there does not exist any non-trivial
generic Alexander polynomials of degree $n$ with
$n\not \equiv 0$ mod 6, for example, this is the case for
 cubic, quartic and quintic curves,
whose singularities are
copies of $A_1$  or  $A_2$. However even though
{\em there does  exist  interesting geometry on these curves}.
We will show by examples that certain tangential Alexander polynomial 
gives non-trivial Alexander invariants and we will
 give an explanation from viewpoint of
non-reduced degeneration in \S 5.
 
Another weakness of generic Alexander polynomials  is  for reducible curves.
Let $C_1$ and $C_2$ be curves which intersect transversely each other.
We  take a  line  $L$  so that
$L$
does not contain any points of $C_1\cap C_2$. 
Note that $L$ need not be generic for $C_1\cup C_2$.
 Theorem \ref{Oka-Sakamoto} says that
\[
 \pi_1(\BC_L^2-C_1\cup C_2)\cong \pi_1(\BC_L^2-C_1)\,\times\,  \pi_1(\BC_L^2-C_2)
\]
However  the Alexander polynomial $\De_{C_1\cup C_2}(t;L)$ loses these
informations. In fact, we have
\begin{Theorem}{\rm  (\cite{OkaSurvey})}\label{Product2} Assume that $C_1$ and $C_2$
 intersect transversely and let $C=C_1\cup C_2$. Let $L$ be a line
such that $L\cap C_1\cap C_2=\emptyset$.
Then $\De_C(t;L)=(t-1)^{r-1}$ where $r$ is the number of irreducible components.
\end{Theorem}
For further information about Alexander polynomials, we refer to
\cite{Degtyarev-alex,Libgober1,Libgober2,Loeser-Vaquie,Randell}
\section{Dual stratification and tangential fundamental groups.}
\subsection{Dual stratification of curves}
Let $\Si$ be a finite set of topological equivalent class of curve
singularities
and let $\cM(\Si,d)$ be {\em the configuration space} of  plane curves of
degree $d$ with a fixed singularity configuration $\Si$. 
Take two curves $C,\,C'\in \cM(\Si,d) $ in the same connected 
component and two smooth
 points $P\in C$ and $Q\in C'$.
We consider their tangent lines
$L=T_PC,\,L'=T_QC$.
Though the topology of 
$(\BP^2,C)$ and $(\BP^2,C')$ are topologically equivalent,
 this may not the case for 
$(\BP^2,C\cup L)$ and $(\BP^2,C'\cup L')$. To analyze this,
we introduce the {\em dual stratification 
$\cS(C)$ for $C\in\cS(\cM(\Si,d))$ and  $\cS(\cM(\Si,d))$
 of $\cM(\Si,d)$}
as follows.

 Let $\check \BP^2$ be the dual projective space. Recall that 
a point $\al=(\al_1:\al_2:\al_3)\in \check \BP^2$ (resp. a point $P=(p_1:p_2:p_3)\in \BP^2$) can be
considered as a line $L_\al=\{\al_1\,X+\al_2\,Y+\al_3\,Z=0\}$ in $\BP^2$ 
(resp. a line $L_P=\{p_1\,U+p_2\,V+p_3\,W=0\}$ in $\check \BP^2$).
First take $C\in\cM(\Si,d) $. Let $\Si(C)=\{P_1,\dots, P_k\}$ be the singular
 points
of $C$.  Let $\cP(d)$ be the set of partition of the integer $d$.
We consider the mappings
$ \psi: \check C \to \cP(d)$ and
$\check \psi: C\to \cP(\check d)$
defined  as follows.
Let $\al\in \check C$ (resp. $P\in C$)
and let $L_\al\cap C=\{R_1,\dots, R_\nu\}$ (resp. 
$L_P\cap \check C=\{S_1,\dots, S_\mu\}$). 
We define $\psi(\al)=\{I(C,L_\al;R_i),\, i=1,\dots, \nu\}$
where $I(C,L_\al;R_i)$ is the local intersection multiplicity.
Respectively we define $\check\psi (P)=\{I(\check C,L_P;S_j);j=1,\dots, \mu\}$.
Note that for a generic line $\al\in \check C$,
$L_\al$ is a simple tangent line and therefore $\psi(\al)=\{2,1,\dots,
1\}$.
For a generic flex point $P$, the tangent line $L=T_P C$ 
gives the partition  $\psi(L)=\{3,1,\dots,1\}$.
A line $\al\in \check C$ is called an {\em multi-tangent} line
if $\psi(\al)$ has at least two members $\ge 2$.  
A simple bi-tangent line is a typical such line
which is simply tangent at two smooth points. A smooth point $P\in C$
is called {\em  tangentially generic} if it is smooth and the tangent
line $T_PC$ 
gives the partition $\{2,1,\dots,1\}$.
 Recall that the Gauss map associated with $C$, denoted as $G_C: C\to \check C$, 
is defined by 
$G_C(P)=T_P C$.
Let $\Si^{ntg}(C)=\{P_{k+1},\dots, P_{k+t}\}$ be smooth points which are
 not tangentially generic and put
$\widetilde {\Si}(C)=\Si(C)\cup  \Si^{ntg}(C)=\{P_1,\dots, P_{k+t}\}$.
The dual stratification $\cS(C)$ of $C$ is  by definition,
$\cS:=\{C-\widetilde {\Si}(C), \widetilde {\Si}(C)\}$. Thus
if $C$ is irreducible,  $\cS(C)$  has 
one open dense  stratum made of tangentially generic points
 and $k+t$ starata made of isolated points.
\subsection{Dual stratification of the configuration space $\cM(\Si,d)$}
Now we consider the dual stratification of $\cM(\Si,d)$.
To distinguish  a point in $\cM(\Si,d)$ and the corresponding  curve, we 
denote points in $\cM(\Si,d)$ by  $\alpha\in \cM(\Si,d)$
and  the corresponding curve by $C_\alpha$.
 The  configuration of the  singularities of the dual curve
$\Si(\check{C_\al})$   is not unique for $\al\in \cM(\Si,d)$ but 
it has only finite possible  types,  say $\Si_1^*,\dots, \Si_{\ell}^*$
when we fix the configuration space  $\cM(\Si,d)$.
We consider the partition of the configuration space by the 
following sets:
\[
 \{\al\in \cM(\Si,d); \left(\Si(\check{C_\al}),\, \cS(C_\al),\,
 \cS(\check{C_\al}), \psi_\al,\check \psi_\al\right)\,\,\text{are constant}\}
\]
The {\em dual
stratification} $\cS(\cM(\Si,d))$
is defined by the strata which are  the connected components
of  these  partitions.
Thus for a stratum $M\in \cS(\cM(\Si,d))$,  each $C_\al$ and $\check{C_\al},\, \al\in M$
have constant dual stratifications.


For a stratum $M\in  \cS(\cM(\Si,d))$, we  can associate a family of
plane curves $C_\al,\,\al\in M$  such that 
the dual family  of curves $\check{C}_\al,\, \al\in M$
is a family in $\cM(\Si_j^*,\check{d})$ for some $j$. 
Observe  that
any $\al\in M$, the dual stratification $\cS(C_\al)$ and $\cS(\check
C_\al)$
are constant for
$\al$
by definition.  Thus for two $\al,\be\in M$,
$C_\al,\,C_\be$ are homeomorphic as a stratified sets.
More precisely we have
\begin{Proposition} Take $\al_0\in M$ and take 
 a  point
  $L_{\al_0}\in \check{C}_{\al_0}$.
This induces a continuous family of lines
 $L_\al\in \check{C}_\al$
such that $\psi_{\al}(L_\al)$ is constant. Then
the topology of the affine pair $(\BC_{L_{\al}}^2,C)$   does not depend on $\al\in M$.
In particular  the  fundamental group
$\pi_1(\BC_{L_{\al}}^2-C)$ does not depend on $\al\in M$.
\end{Proposition}
{\em Proof}. Recall that the local topology of
$C_\al\cup L_\al$ at an intersection point $P$
is determined by the local Milnor number  $\mu(C_\al\cup L_\al,P)$
and this is determined by $\mu(C_\al,P)$
and  the local  intersection multiplicity
$I(C_\al,L_\al;P)$. The definition of the dual stratification 
of $\cS(\cM(\Si,d))$ guarantees the $\mu$-constancy of the family
of plane curves
$C_\al\cup L_{\al},\, \al\in M$ of degree $d+1$. \qed


\subsection{Tangential fundamental group and Tangential Alexander
  polynomial}
For a line $L\in \check C$, 
we call $\pi_1(\BC_L^2-C)$ {\em the tangential fundamental group} and 
$\De_C(t;L)$  {\em the tangential Alexander polynomial}.
If   
 $L=T_PC$ for some simple point
$P\in C$, we also use the notation 
$\De_C(t;P)$ for $\De_C(t;T_PC)$.
We also define $k$-fold Alexander polynomial
$\De_C(t;P_1,\dots,P_k)$ by
$\De_{C\cup L_1\cup \cdots\cup L_{k-1}}(t;P_k)$ with $L_j=T_{P_j}C$.
It is easy to observe that $\pi_1(\BC_L^2-C)$ and $\De_C(t;P)$ 
are constant on the open (dense if $C$ is irreducible)
strata
of the dual stratification $\cS(C)$.
However  in general it may give a different polynomial for singular
lines $L\in \check C$ (they are the images of isolated strata of
$\cS(C)$ by the Gauss map).
We will see some examples later.
Thus the tangential Alexander polynomials altogether contain more
 geometrical informations than the generic Alexander polynomials.
The main purpose of this paper is  to investigate the property of
the tangential Alexander polynomials. Note that if $C$ is
irreducible,
there is only one choice of $\theta$ (up to $\pm$) but there are  many
choices for $L$, even for irreducible $C$.

\subsection{Alexander spectrum}
We also
  consider  the set of tangential Alexander polynomials
\[
 \tAS\,(C):=\{\De_C(t;P);P\in C\}
\] 
and we call $ \tAS\,(C)$
  {\em the  tangential Alexander spectrum of $C$}.
There exist at most finite polynomials in the spectrum.
In fact, it is bounded by the number of strata of $\cS(C)$.

We can also  define the {\em $k$-fold tangential Alexander spectrum of $C$}
by
\[
 \tAS^{(k)}(C):=\{\De_{C}(t;P_1,\dots,P_k); P_j\, \in C\}
\]
It often happens that even when the Alexander spectrum $\tAS(C)$
is
trivial, 2-fold Alexander spectrum  $\tAS^{(2)}(C)$ (or higher one) is not trivial.

\subsection{Example} \label{example-quartic}

We consider $\cM(2A_2+A_1,4)$ and $\cM(E_6,4)$. By class formula,
the dual curve $\check C$ 
of a generic member $C$ of $\cM(2A_2+A_1,4)$ or $\cM(E_6,4)$
 is a quartic with $2A_2+A_1$ in both cases. ($C$ has generically 2 flex points.)
In both configuration spaces
$\cM(2A_2+A_1,4)$ and $\cM(E_6,4)$, there are  strata   which
correspond to {\em degenerated} members, namely curves with
 one flex of order 2. This implies that  the dual curve has 
an $E_6$ singularity. For these configuration spaces, there  is a 
beautiful work by C.T.C.Wall
\cite{Wall-quartic}.
Consider the subsets:
\begin{eqnarray*}
& M_1:=\{C\in \cM(2A_2+A_1,4); \Si(\check C)=\{2A_2+A_1\}\},\quad\\
&M_2:=\{C\in \cM(2A_2+A_1,4); \Si(\check C)=\{E_6\}\}\\
& N_1:=\{C\in \cM(E_6,4); \Si(\check C)=\{2A_2+A_1\}\},\quad\\
&N_2:=\{C\in \cM(E_6,4); \Si(\check C)=\{E_6\}\}
\end{eqnarray*}
We can easily see that 
$\{M_1,\,M_2\},\,\{N_1,\,N_2\}$ are respective dual stratifications
of the configuration spaces $\cM(2A_2+A_1,4)$ and $\cM(E_6,4)$.
We observe that under the Gauss map, 

--$M_1$ and $N_2$ are self-dual and 

--$M_2$ and $N_1$ are dual each other.

\noindent
We observe also that $M_2 \subset \partial M_1$
and $N_2\subset \partial N_1$. 

\vspace{.5cm}\noindent
(1-1) We consider quartic $C_1\in M_1$ with $\Si(C_1)=2A_2+A_1$ with 
two flexes. By the class formula, such a curve has a bi-tangent line. As
an example, we take:
\begin{eqnarray}
C_1:\,{\frac {17}{4}}\,{y}^{4}+8\,{y}^{3}+1/4+7/2\,{y}^{2}-7/2\,{y}^{2}{x}^{
2}+1/4\,{x}^{4}-1/2\,{x}^{2}=0
\end{eqnarray}
$C_1$ has two cusps at $P_1=(1,0),\,P_2=(-1,0)$ and one  $A_1$ at
$(0,-1)$. Two flexes are at
$Q_1=(10\sqrt{10},9),\, Q_2=(-10\sqrt{10},9)$. 
We have a bi-tangent line $y=1$ which are tangent at
$B_1=(2\sqrt 2,1),\,Q_2=(-2\sqrt 2,1)$. For the dual stratification
$\cS(C_1)$, we  have to take two more points
$S_1=(25/7,8/7),\,S_2=(-25/7,8/7)$ whose tangent lines pass through
$P_2$ and $P_1$ respectively. Using Zariski-van Kampen pencil
method,
we can compute $\pi_1(\BC_L^2-C_1)$ as
\[
 \pi_1(\BC_L^2-C_1)=\langle
\xi_1,\xi_2,\xi_3\,|\,\, \xi_1\xi_2\xi_1=\xi_2\xi_1\xi_2,\,
\xi_2\xi_3\xi_2=\xi_3\xi_2\xi_3,\,\xi_1\xi_3=\xi_3\xi_1,\,
\xi_1\xi_2\xi_3=\xi_3\xi_2\xi_1\rangle
\]
where $L=\{y=1\}$. This gives $\De_{C_1}(t;L)=t^2-t+1$ by Fox calculus.
Other tangent lines give the trivial Alexander polynomial.
We leave the proof of this assertion as an exercise.
Thus $\tAS(C_1)=\{1,t^2-t+1\}$.

\noindent
(1-2) We consider the following quartic $C_2\in M_2$ with $\Si(C_2)=2A_2+A_1$ with a
 degenerated flex
 of order 2 at infinity $L=\{Z=0\}$:
\[C_2:\quad
1-12\,y+36\,{y}^{2}-32\,{y}^{3}-2\,{x}^{2}+12\,{x}^{2}y+{x}^{4}=0
\]
\[
 \pi_1(\BC_L^2-C_2)=\langle
\xi_1,\xi_2,\xi_3| \xi_1\xi_3\xi_1=\xi_3\xi_1\xi_3,\,  \xi_3\xi_2\xi_3=\xi_2\xi_3\xi_2,\,\xi_2\xi_1=\xi_1\xi_2\rangle
\]
We can also see that $\De_{C_2}(t;L)=t^2-t+1$.
Note that two $A_2$ singularities are at $(\pm 1,0)$ and one $A_1$ is at $(0,1/2)$.
In the dual stratification $\cS(C_2)$, there are two more 'singular'
points
$R_1=(9,8)$ and $R_2=(-9,8)$ whose tangent line pass through the cusps.
However these tangent lines give trivial tangential Alexander polynomials.

\noindent
(2-1) Consider a quartic $C_3\in N_1$, defined by
$ y^3+x^4-x^2y^2=0$
which has one $E_6$-singularity at $O$. Two flex points are at
$(\pm 6\sqrt 6/5,36/5)$.
The bi-tangent line is given by $y=4$. By an easy computation, we observe
that
$L=\{y=4\}$ gives $\pi_1(\BC_L^2-C_3)\cong B_3$ and 
 $\De_{C_3}(t;L)=t^2-t+1$. Other tangent Alexander
spectra are trivial. So $\tAS(C_3)=\{1,t^2-t+1\}$.
The dual stratification has 
$5$ isolated points.

\noindent
(2-2) Consider the following  quartic $C_4\in N_2$,
$f(x,y)=\, y^{3} + x^{4}=0$, with $E_6$ and one flex
of order 2 at  $P=(0,1,0)$. 
Take the flex line $Z=0$ as $L$. Then
\[
  \pi_1(\BC_L^2-C)\cong G(3,4)=
\langle \xi_0,\xi_1,\xi_2|\, \xi_0=\omega\xi_1\omega\inv,\,
\xi_1=\omega\xi_2\omega\inv\rangle
\]
and $\De_{C_4}(t;L)=(t^2-t+1)(t^4-t^2+1)$
where $\omega=\xi_2\xi_1\xi_0$. 

\noindent
(3) Let $C:\,f(x,y) = - {\displaystyle \frac {1}{2}} \,y^{4} - {\displaystyle 
\frac {1}{2}}  - 3\,x^{2}\,y^{2} + y^{2} + {\displaystyle \frac {
3}{2}} \,x^{4} - 4\,x^{3} + 3\,x^{2}=0$ be a 3 cuspidal quartic.
 As is well-known \cite{Za1,OkaSymmetric}, the generic affine
 fundamental group is a finite group of order 12, with presentation
\[
\pi_1(\BC^2-C)=\langle\xi,\zeta|\xi\,\zeta\,\xi=\zeta\,\xi\, \zeta,\,\xi^2=\zeta^2\rangle.
\]
Though the fundamental group is not abelian,
 the generic Alexander polynomial is trivial.
For $L=\{y=0\}$ (this is the tangent cone of a cusp and $L$
corresponds to a
 flex of $\check C$),
\[
 \pi_1(\BC_L^2-C)=\langle\xi,\zeta| \xi\,\zeta\,\xi=\zeta\,\xi\,
 \zeta\rangle=B_3.
\]
By the class formula, the dual curve $\check C$ is a cubic curve with a node.
$\cS(C)$ has 3 singular points  from $3A_2$ and two 'singular points'
from   the bi-tangent line.
We can also see that $\pi_1(\BC_L^2-C)\cong B_3$  
for  an arbitrary tangent line $L$ except the bitangent line $L_b$.
The bi-tangent line is given by $x=2/3$ in this example.
By an easy computation, we see that
\[
 \pi_1(\BC_{L_b}^2-C)\cong\langle \xi_0,\xi_1,\xi_2,\zeta\,|\,
\xi_0\xi_1\xi_0=\xi_1\xi_0\xi_1,\,\xi_1\xi_2\xi_1=\xi_2\xi_1\xi_2,\,
\xi_2\zeta\xi_2=\zeta\xi_2\zeta,\,\zeta=\xi_1\inv\xi_0\xi_1\rangle
\]
and
we have
\begin{Proposition}\label{Exceptional} For the bitangent line $L_b$,
  $\De_{C}(t;L_b)=(t^2-t+1)^2$. For any other tangent line $L$,
$\De_C(t;L)=t^2-t+1$.
In particular, this implies that
$\tAS(C)=\{t^2-t+1,\,(t^2-t+1)^2\}$.
\end{Proposition}

\subsubsection{Further example} In the above examples of quartics, the
geometry of  $C\cup T_PC$ does not change for flexes of the same order.
However this is not the case in general. 

\noindent
Consider a fixed mark  point $P\in C$.
 We call $(C,P)$ {\em a  curve with a marked point $P$}.
Two curves with marked points $(C,P)$ and $(C',P')$ are called 
 {\em a marked Zariski pair} if
$\{C\cup T_PC, C'\cup T_{P'}C'\}$ is a Zariski pair.
For  further information about Zariski pairs, see \cite{Artal,MR99d:14019,MR2001a:14024,OkaSurvey,ZariskiPairsI}.
In \cite{ZariskiPairsII}, we have shown that
for any  quintic $B_5$ with 
configuration in the next  list,
there exist
 two different flex points
$P,\,P'$ such that $(B_5,P)$ and $(B_5,P')$ are marked Zariski pairs.
\begin{eqnarray*}
(\sharp)\qquad 
\begin{cases}
& 4A_2\,, \,4A_2+A_1\,, \,A_5+2A_2\,, \,A_5+2A_2+A_1\,, \,E_6+2A_2 \\
&E_6+A_5, \,2A_5,\,
A_8+A_2\,, \,A_{8}+A_2+A_1, \,A_{11}
\end{cases}
\end{eqnarray*} 
In fact their generic Alexander polynomials
are given as 
$t^2-t+1,\, 1$ respectively.
This implies  that the among flexes of these quintics, there are two
classes of different topological nature: one class which does not contribute the tangential Alexander
spectrum
and the other which contributes by $(t^2-t+1)$.
We give one example.
The following quintic $B_5:\,f(x,y)=0$ has $A_{11}$  singularity
at the origin and  9
 flex points. Among them, the flex at $P=(0,1)$ is different from others
(a flex of torus type).
In fact, $B_5\cup T_PB_5$ is a sextic of torus type \cite{ZariskiPairsII}.
All other flex points gives trivial tangential Alexander polynomial.
\begin{multline*}
f(x,y)=-{\frac {33}{64}}\,{y}^{5}+\left( {\frac {7}{8}}\,x+{\frac {129}{64}}
 \right) {y}^{4}+ \left( -\frac 54\,{x}^{2}-{\frac {15}{8}}\,x-5/2 \right) 
{y}^{3}+\\
 \left( {\frac {15}{8}}\,{x}^{3}
+{\frac {13}{4}}\,{x}^{2}+x+1
\right) {y}^{2}
+ \left( -\frac 34\,{x}^{4}-2\,{x}^{3}-2\,{x}^{2} \right) y
+{x}^{5}+{x}^{4}
\end{multline*}
\begin{figure}[htb]
\setlength{\unitlength}{1bp}
\begin{picture}(600,150)(-100,0)
\put(20,130){\special{epsfile=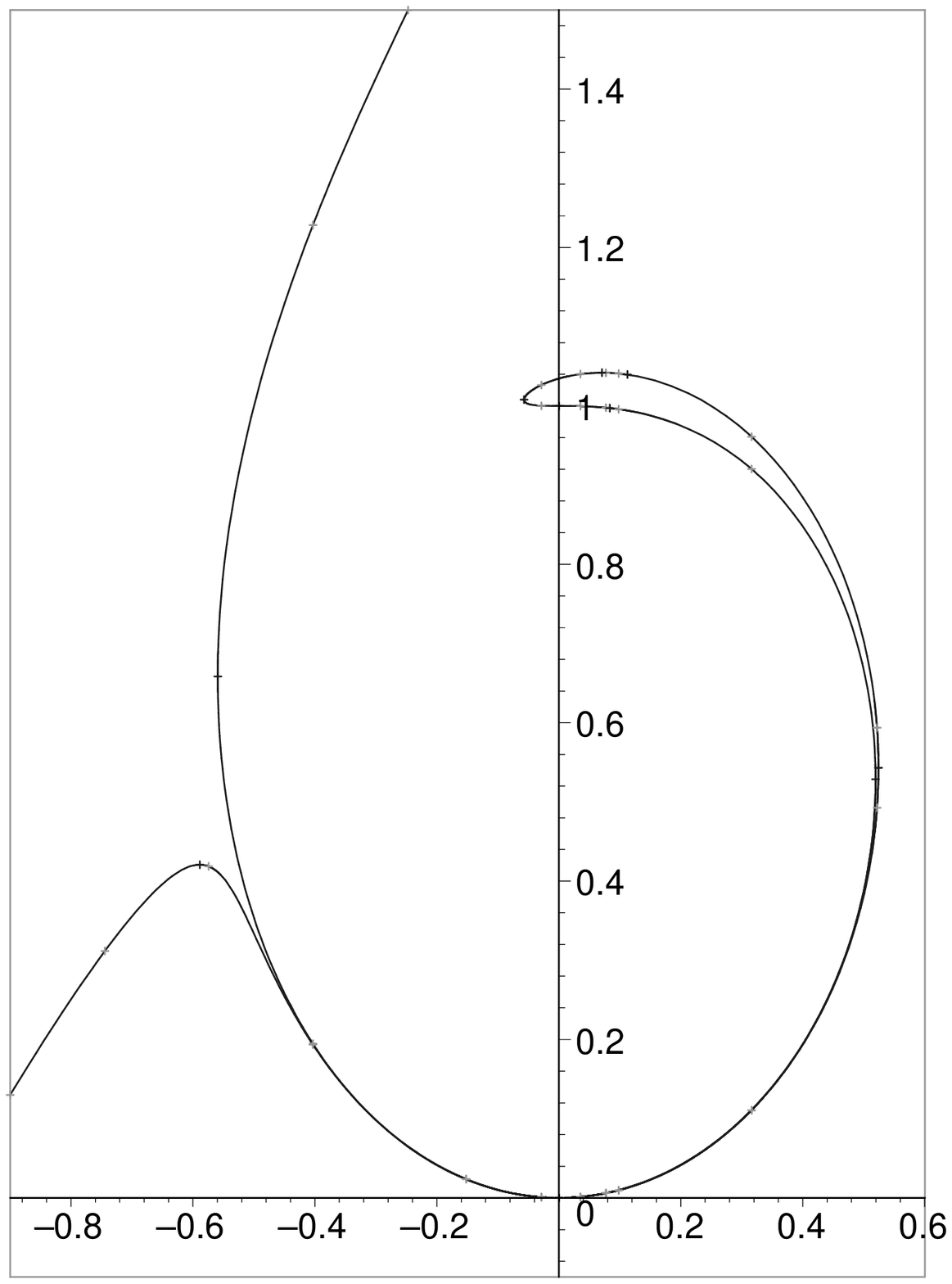 vscale=0.3 hscale=0.4}}
\end{picture}
\vspace{.6cm}
\caption{Quintic with $A_{11}$}
\end{figure}
\noindent

\subsection{$\theta$-Alexander polynomials}
To cover the weakness of Alexander polynomials
for irreducible curves, we have
proposed $\theta$-Alexander polynomials in \cite{OkaSurvey}.
First recall that  the {\em  radical}  $\sqrt{q(t)}$ of a polynomial 
$q(t)=\prod_{i=1}^\nu (t-\xi_i)^{\mu_i}$ is defined by
  $\sqrt{q(t)}:=\prod_{i=1}^\nu (t-\xi_i)$.
Here $\mu_i\ge 1,\forall i$. 
The following theorem shows the importance of $\theta$-Alexander polynomial.
\begin{Theorem}(\cite{OkaSurvey})\label{Non-trivial}
Assume that $C,\,C'$   be reduced curves and
we assume further $C'$ is irreducible. 
For a given integer $n$, suppose that the surjective homomorphism
$ \phi_n:\pi_1(\BC_L^2-C\cup C')$\newline
$\to \BZ$
which has weight 1 on each component of $C$ and weight $n$ on $C'$.
Then 

\noindent
 (1) $\De_{C\cup C'}(t;L,\phi_n)$ is divisible by
 $\gcd({\De_C(t;L)\times (t-1)},(t^n-1))$
for $n\ne 0$.

\noindent
Suppose further that $C\pitchfork C'$
and $\pi_1(\BC_L^2-C')\cong \BZ$. Then

\noindent
(2) $\De_C(t;L)\times (t-1)$ is divisible by $\De_{C\cup C'}(t;L,\phi_n)$.

\noindent
(3)  $\sqrt{\De_{C\cup C'}(t;L,\phi_n)}=\gcd(\De_C(t;L)\times (t-1),(t^n-1))$.

\noindent
(4)
In particular, if $\gcd({\De_C(t;L)\times (t-1)},(t^n-1))=\De_C(t;L)$,
we have
$$\De_{C\cup C'}(t;L,\phi_n)=\De_C(t;L)\times (t-1).$$
\end{Theorem}
{\bf Proof.}
The proof goes exactly as in \cite{OkaSurvey}.
Consider the canonical surjective homomorphism:
$h:\pi_1(\BC_L^2-C\cup C')\to \pi_1(\BC_L^2-C)\times \BZ$.
Consider the presentation.
\[
 \pi_1(\BC_L^2-C)=\langle g_1,\dots,g_s\,|\,R_1,\dots,
 R_k\rangle,\,\,s\le
\degree\,C
\]
The homomorphism $\phi_n$ factors as $\phi_n=\psi\circ h$ where
$\psi$ is the surjective homomorphism
\newline
$ \psi:\,\,\pi_1(\BC_L^2-C)\times H_1(\BC_L^2-C')\cong
\pi_1(\BC_L^2-C)\times \BZ\to \BZ$ where the weight of  the second factor
is $n$.
Note that 
$\De_{C\cup C'}(t;L,\phi_n)=\De_{\phi_n}(t)$ in the notation of \S 3.1.
By
 the above factorization, we have
the divisibility:
$\De_\psi| \De_{\phi_n}$. Now the calculation of $\De_\psi$ is done in the
exact same manner as in \cite{OkaSurvey}. We use the presentation:
\[
 \pi_1(\BC_L^2-C)\times \BZ=\langle
g_1,\dots, g_s,\xi\,|\,R_1,\dots, R_k,T_{i},\,
1\le i\le s 
\rangle
\]
where $T_{i}=g_i\xi g_i\inv \xi\inv$.
The key point of the calculation is the following:
\begin{eqnarray*}
\frac{\partial T_{i}}{\partial g_i}=1-g_i\xi g_i\inv\mapsto
1-t^n,\quad
\frac{\partial T_{i}}{\partial \xi}=g_i-g_i \xi g_i\inv \xi\inv\mapsto
t-1
\end{eqnarray*}
Let $M$ be the Alexander matrix of 
$\theta_{sum}:\pi_1(\BC_L^2-C)\to \BZ$
and let $M'$ be the Alexander matrix of 
$\psi:\pi_1(\BC_L^2-C)\times \BZ\to \BZ$. Then
$M'$ is written as
\[
M'=\left( \begin{matrix}
M&\vec 0\\
N_1&N_2\\
\end{matrix}
\right)
\]
where $\vec 0$ is a zero vector and $N_1$
is a $s\times s$-matrix which 
is explicitly given as $(1-t^n)E_{s}$ where $E_{s}$ is the  identity
matrix of rank $s$.
 The vector $N_2$ takes the form
${}^t(t-1,\dots, t-1)$
where $t-1$ is repeated $s$ times. For
any $(s-1)\times (s-1)$--minor $B$ of $M$, let 
$\tilde B$ be the $s\times s$--minor adding $(k+1)$-th row and the last
column.
Then $\det\,\tilde B=\det\, B\times (t-1)$.
Thus we have
 $\De_{C}(t;L)\times (t-1)$
as the common divisor of such  $\det\,\tilde B$'s.
Also we get  $(t^n-1)^{s}$ by taking a minor  from $N_1$.
We observe also  that  any  determinant of a  $s\times s$--minor
which contains at least two rows of $(N_1,N_2)$ is divisible by $(t^n-1)$.
 Thus  we observe two divisibilities: 
\begin{eqnarray*}
\De_{\psi}(t)\,|\,\gcd( \De_{C}(t;L)\times (t-1),\, (t^n-1)^{s}),\quad
\gcd(\De_{C}(t;L)\times (t-1), (t^n-1))\,|\, \De_{\psi}(t)
\end{eqnarray*}
Note that 
$\De_{\psi}(t)\,|\, \De_{C\cup C'}(t;\phi_n,L)$,
by  the usual degeneration argument \cite{OkaSurvey}.  
Thus the first assertion is
immediate from the last divisibility.
Suppose further that $C\pitchfork C'$ 
and $\pi_1(\BC_L^2-C')\cong \BZ$.
Then
$h$ is an isomorphism and therefore $\De_{C\cup C'}(t;\phi_n,L)=\De_{\psi}(t)$.
Thus the assertions (2), (3) follow immediately. The assertion (4) is a
result of (1) and (2).
\qed

\subsection{Relations between the tangential Alexander polynomials and $\theta$-Alexander
 polynomials}

Let $C$ be a plane curve of degree $d$ and let $P\in C$  and let $L=T_PC$.
We consider the tangential Alexander polynomial
$\De_C(t;L)$. Let 
\[
\pi_1(\BC_L^2-C)\cong  \langle g_1,\dots,g_d\,|\, R_1,\dots, R_\ell\rangle
\]
be a presentation of $\pi_1(\BC_L^2-C)$ by generators and relations.
Take a generic line $L_\infty$ for $C\cup L$ and put
$\BC^2=\BP^2-L_\infty$ as usual.
Then by Theorem \ref{Oka-Sakamoto},
we have
\begin{multline*}
 \pi_1(\BC^2-C\cup L)=\pi_1(\BC_L^2-C\cup L_\infty)
\cong\pi_1(\BC_L^2-C)\times\pi_1(\BC_L^2-L_\infty)
\cong\pi_1(\BC_L^2-C)\times\BZ
\end{multline*}
and it has a presentation:
\begin{eqnarray*}
\pi_1(\BC^2-C\cup L)=& \langle g_1,\dots,g_d,h, h_\infty\,|\,
 R_1,\dots,R_\ell, T_1,\dots,
 T_d, S\rangle\\
=&\langle  g_1,\dots,g_d,\,h_\infty\,|\,
R_1,\dots,R_\ell, T_1,\dots,
 T_d
\rangle
\end{eqnarray*}
where $T_j=\, h_\infty\, g_j\, h_\infty\inv\, g_j\inv$
 and 
$S=\, h\, h_\infty\,g_d\dots g_1$.
Now the tangential Alexander polynomial
$\De_C(t;L)$ is associated to the 
surjective  homomorphism
\[
 \theta_{sum}:\pi_1(\BC_L^2-C)\to \BZ=\langle t\rangle,\, g_i\mapsto t
\]
Let $\theta_n$
be the surjective
homomorphism with weight $n$ on $L$
\[
\theta_n:\pi_1(\BC^2-C\cup L)\mapright{} \BZ=\langle t\rangle,\,
g_i\mapsto t,\, h\mapsto t^n
\]
Now taking $g_1,\dots,g_k,\, h_\infty$ as generators
of $\pi_1(\BC^2-C\cup L)=\pi_1(\BC_L^2-C\cup L_\infty)$,
$\theta_n$ corresponds to the homomorphism:
\[
 \eta_{n+d}:\pi_1(\BC_L^2-C\cup L_\infty)\to \BZ,\quad
g_i\mapsto t,\quad  h_\infty\mapsto  t^{-n-d}
\]
The last property is  the result of the observation:
$1=\theta_n(S)=t^{d+n} \theta_n(h_\infty)$.

\noindent
{\bf Notation.} Hereafter we mainly consider the weight like $\theta_n$
which has  weight one  except a line component $L$ in consideration. So we introduce the
following notation which is  easier to be understood:
\[
 \De_{C\cup L^n}(t):=\De_{C\cup L}(t;\theta_n)
\]
The upper index $n$ implies that $L$ has weight $n$.
Using  this notation, we  also write
\newline
 $\De_{C\cup L_{\infty}}(t;,\eta_{n+d},L)\,=\,\De_{C\cup L_\infty^{-n-d}}(t;L)$.

\vspace{.2cm}\noindent
Thus combining the above argument  with Theorem \ref{Non-trivial} we
have shown
the following.
\begin{Theorem}\label{Tangential}
For any integer $n$,
$\De_{C\cup L^n}(t)= \De_{C\cup L_\infty^{-n-d}}(t;L)$ and we have the divisibility:
\[
 \De_{C\cup L_\infty^{-n-d}}(t;L)|\De_C(t;L)\times (t-1)
\quad\text{and}\quad
\gcd (\De_C(t;L)\times (t-1),(t^{n+d}-1))|\De_{C\cup L_\infty^{-n-d}}(t;L)
\]
Furthermore if $\gcd (\De_C(t;L)\times (t-1),(t^{n+d}-1))=\De_C(t;L)$, we have the
 equality:
\[
 \De_C(t;L)\times (t-1)=\De_{C\cup L_\infty^{-n-d}}(t;L).
\]
\end{Theorem}

\subsubsection{Examples} 
Let $C$ be a quartic with either $2A_2+A_1$ or $E_6$ and one flex of
 order 2. Let $L$ be the flex tangent line. Then we have shown that 
$\De_C(t;L)=t^2-t+1$ and  $(t^2-t+1)(t^4-t^2+1)$
respectively. 
We can 
 compute  generic 
Alexander polynomials $\De_{C\cup L^2}(t)$ of $C\cup L$ as follows. Take a
 generic line $L_\infty$.

\vspace{.2cm}
\noindent
(1) $C$ is a quartic with  $2A_2+A_1$ and $L$  is the flex tangent line.
As $\De_{C}(t;L)=t^2-t+1$,  we take
weight 
$n=2$  on $L$ 
and by Theorem \ref{Non-trivial}, we get
\[
 \De_{C\cup L_\infty^{-6}}(t;L)=
\De_{C\cup L^2}(t)=(t^2-t+1)(t-1)
\]
 Let $C_t$ be a family of quartics with
$2A_2+A_1$ with two flex points for $t\ne 0$ and $C_0=C$.
Let $L_1,\,L_2$ be two tangent lines at the flex points of $C_t$.
Then  $C_t+L_1+L_2\to  C+2L$.
Thus the weight 2 on $L$ is canonical. See \S 5.7.

\noindent
(2) $C$ is a quartic with  $E_6$ and $L$  is the flex tangent line at a
flex of order 2.
As $\De_C(t;L)=(t^2-t+1)(t^4-t^2+1)$, 
$\gcd(\sqrt{\De_C(t;L)},(t^{12}-1))=\De_C(t;L)$.
Thus we take $n=8$.
\[
 \De_{C\cup L^8}(t)=\De_{C\cup L_\infty^{-12}}(t;L)=
(t-1)(t^2-t+1)(t^4-t^2+1)
\]
This can be interpreted as $Y^3Z+X^4=0$ is  a line degeneration of (3,4)-torus
curve 
of degree 12 as 
$(Y^3Z+X^4)Z^8=(YZ^3)^3+(XZ^2)^3$. See \S 5.4.
Note also that 
\[
\De_{C\cup L^2}(t)= \De_{C\cup L_\infty^{-6}}(t;L)=(t-1)(t^2-t+1).
\]
We can also interpret this equality
 as a result of a line degeneration of (2,3)-sextics
of torus type as
 $(Y^3Z+X^4)Z^2=(YZ)^3+(X^2Z)^2$.

\section{Degeneration into non-reduced curves with a multiple line}
In this section, we  study  an analytic 
family of
curves $C_t,\,t\in \De$ such that  $C_0$ is not reduced but it
has a line component with multiplicity.

\subsection{Admissible polydisk}
Consider a reduced  curve
$C\subset \BP^2$ which is defined by a polynomial $f(x,y)=0$ in the
affine space $\BC_L^2:=\BP^2-L$ where 
$L=\{Z=0\}$. (We do not assume the genericity of
the line $L$.)
We assume that $f(x,y)$ is a   polynomial in $y$ of degree $n$.
The base point of the pencil $\{L_\eta,\eta\in \BC\}$
where $L_\eta:=\{X-\eta Z=0\}$
is given by $B=(0,1,0)$ in the homogeneous
coordinates. Note that $n<d$ if and only if $B\in C$.
We say the pencil $\{x=\eta,\,\eta\in \BC\}$ is {\em base point
admissible}
(respectively {\em base point non-admissible}) if $n=d$  (resp. $n<d$).
Recall that $L_\eta$ is a singular line for $C$ if $L_\eta\cap C\cap \BC_L^2$
contains some non-transverse intersection
 point. For the case $n<d$, we also call
$L_\eta$
singular if the number of the points $L_\eta\cap C\cap \BC_L^2$
counted with multiplicity is strictly
less than $n$. In this case, we say $L_\eta$ {\em a singular line with
disappeared points at infinity}.
Using Zariski-van Kampen  pencil method with respect to
the pencil lines $x=\eta,\,\eta\in \BC$, we get a
presentation
\[
 \pi_1(\BC_L^2-C)=\langle
g_1,\dots, g_n\,|\, R_1,\dots, R_m\rangle
\]
We consider the  polydisk $ \Delta_{\al,\be}:=\De_\al \times \De_\be$
where $\De_\al :=\{x\,|\, |x|\le \alpha\}$ and
$\De_\be :=\{y\,|\,|y|\le \beta\}$
 and we consider the following conditions.

\begin{enumerate}
\item For
any $\eta\in \De_\al $, 
 $L_\eta\cap C\subset \{\eta\}\times\De_{\be/2} $
and $L_\eta$ does not have  any
disappeared points at infinity and 
\item For any $\eta\in \partial \De_\al $, $L_\eta$ is not a singular
line.

\item For any  singular line  $L_\eta$, $\eta\in \De_\al.$
\end{enumerate}
We say that the polydisk $\De_{\al,\be}$ is
{\em admissible for $C$ with respect to $L$}
 if it satisfies (1) and (2).
Furthermore, we say that the polydisk $\De_{\al,\be}$ is
{\em topologically presenting for $C$ with respect to $L$}
 if it satisfies (1), (2) and (3).
Observe that if $\De_{\al,\be}$ is 
 topologically presenting for $C$ with respect to $L$,
 the inclusion
$ \De_{\alpha,\beta}-C\cap  \De_{\alpha,\beta}\hookrightarrow
\BC_L^2-C$
is a homotopy equivalence.

\subsection{Two non-reduced degenerations}In this section, we focus 
the following two types of non-reduced degenerations.

\noindent
{\bf  Type 1: Line Degenerations}. $\{C_t,\,t\in \De\}$ is an analytic family of 
 irreducible curves
of degree $d$ and it degenerate into
 $ C_0=D_0+jL_\infty,\,j\ge 0$ where
$D_0$ is an irreducible curve of degree $d-j$
and $L$ is a line. We assume also 
\newline\indent
($\sharp$) there is a 
point $Q\in L-L\cap D_0$
such that 
$Q\in C_t,\, \forall t\ne 0$
and the multiplicity of
\newline\indent
 $(C_t,Q)$ is $j$. 
\newline
We call such a degeneration
{\em a line degeneration of order $j$}.
$L$ and $Q$ are called the {\em limit line} and the {\em base point}
of the degeneration respectively. Here $j$ is a non-negative
integer.
(We are mainly interested in the
case
$j\ge 2$.)
The condition $(\sharp)$ can be weakened as 
\newline\indent
$(\sharp')$
 there is an analytic family of 
points $Q_t\in L\cap C_t$ such that $Q_0\in L-L\cap D_0$
\newline\indent
and the multiplicity of $(C_t,Q_t)$ is $j$.
\newline\noindent
In fact, under ($\sharp'$),
we may  assume that $L=\{Z=0\}$ and  $Q_t=(\alpha(t),1,0)$.
 Then taking
a linear change of coordinates 
$(x,y)\mapsto(x,y+\alpha(t) x)$, we can assume  that  $Q_t\equiv (0,1,0)$.
We recall that O. Zariski  observed that 3 cuspidal quartic is a
 non-reduced
line degeneration of order 2 from a family of sextics of torus type \cite{Za1}.

\noindent
{\bf Type 2: Flex Degenerations}.  First we have a family of reduced  curves
$\{C_t,\,t\in \Delta\}$ of the same degree. $C_0$ can be a degeneration
if  $\mu(C_0)>\mu(C_t),\,t\ne 0$ or in the same configuration space
if $\mu(C_t)=\mu(C_0)$ 
but in this case, $C_0$ is in a different stratum of the dual stratification.
On $C_t$, we are given 
flex lines $L_1(t),\dots, L_k(t),\,k\ge 2$
such that the family $L_i(t)$
is an analytic family with $L_i(0)=L$ for $i=1,\dots, k$.
We associate a non-reduced  degeneration \newline
$C_t+L_1(t)+\dots+L_k(t)\to
C_0+kL$
and we call this family
{\em a flex degeneration.}

\subsection{Surjectivity Theorem for line degenerations}
Assume first  that we have an analytic family of curves
$\{C_t,\,|t|\le 1\}$ such that  $C_t$ is an irreducible curve
of degree $d$ and it degenerates into
 $ C_0=D_0+jL$ where
$D_0$ is an irreducible curve of degree $d-j$
and $L$ is a line. We assume that $\{C_t,\,t\ne 0\}$ has  
 $Q= (0,1,0)$  as the base point and the multiplicity of $C_t$ at $Q$ is constantly
equal to $j$. Let 
$F(X,Y,Z,t)=0$ be the defining homogeneous polynomial of degree $d$.
We assume that $L=\{Z=0\}$.
By the assumption, we have
$F(X,Y,Z,0)=Z^jG(X,Y,Z)$ where $\degree \, G(X,Y,Z)=d-j$.
Put $f(x,y,t):=F(x,y,1,t)$. This is  the affine equation of $C_t$.
Note that 
\begin{enumerate}
\item $\degree_{\{x,y\}} f(x,y,t)=d$ for $t\ne 0$
\item $\degree_{\{x,y\}} f(x,y,0)=d-j$ and  $\degree_y\, f(x,y,t)=d-j$ for
      any $ t$.
\end{enumerate}
The second assertion follows from   B\'ezout theorem and
 the assumption that  $C_t$ has
multiplicity $j$ at $Q$.
\newtheorem{SurTh}[Theorem]{Theorem (Surjectivity)}
\begin{SurTh}\label{surjectivity}
Under the above assumption, there is a canonical surjection
\[
 \phi:\,\, \pi_1(\BC_L^2-D_0)\,\to\, \pi_1(\BC_L^2-C_\tau),\quad
 \tau\ne 0,\,\text{sufficiently small}
\]
\end{SurTh}
\noindent
{\em Proof.}
Note that the linear system $L_\eta=\{x=\eta z\},\,\eta\in \BC$ has the base
point $Q=(0,1,0)$.
Suppose that we have chosen a topologically presenting
polydisk $\De_{\al,\be}$ for $D_0$.
Let $f(x,y,\tau)$ be the defining affine polynomial for $C_\tau$.
As the effect of the non-reducedness disappears
 when we put $z=1$ in $F(X,Y,Z,t)$, $f(x,y,t)$ is a
analytic family.
Write $f(x,y,t) $ as
\[
 f(x,y,t)=a_{d-j}(x,t)y^{d-j}+\cdots+ a_0(x,t)
\]
 By continuity, we can assume that
this polydisk
 is also admissible for $C_\tau$ for some $\de$ and  $|\tau|\le \de$. 
Let $\eta_1,\dots, \eta_\nu$ be the parameters corresponding to the
singular lines for $D_0$. We take a small positive nuber $\eps$ so that
the disks $\De_\eps(\eta_i):=\{\eta||\eta-\eta_i|\le \eps\},\,i=1,\dots,
\nu$
are disjoint each other and $\De_\eps(\eta_i)\subset \De_\al$.
 By the assumption, we have $|\eta_j|<\al$
for each $j$. We choose a generic pencil line $L_{\eta_0}$
and generators $g_1,\dots, g_d$ in this pencil so that we have a
presentation:
\begin{eqnarray}\label{relation1}
\pi_1(\BC_L^2-D_0)\cong
\langle g_1,\dots, g_d\,\,|\,\, R_1,\dots, R_m\rangle
\end{eqnarray}
For sufficiently small $\tau$, the original singular pencil $L_{\eta_i}$
for
$D_0$
 may splits into
several singular lines for the curves $C_\tau,\, \tau\ne 0$ but  they are inside $\De_\eps(\eta_i)$.
Put the corresponding parameters $\eta_{i,1},\dots,\eta_{i,\nu_i}$.
Note that the monodromy
relations
 around $L_{\eta_i}$ for $D_0$ is nothing but the product of the monodromy
 relations around  $L_{\eta_{i,s}}$ for $ s=1,\dots, \nu_i$ under a suitable
ordering. This is immediate from the topological stability 
of the pencil restricted on the circle $\partial \De_\eps(\eta_i)$.
Note also that $C_\tau$ may have a singular line $L_\eta$ such that
$|\eta|\to \infty$ when $\tau\to 0$.
Anyway we can get a presentation
by adding several more relations
$R_{m+1},\dots, R_{m+n}$ to (2) (using the same generators $g_1,\dots, g_d$):
\[
 \pi_1(\BC_L^2-C_\tau)=
\langle g_1,\dots, g_d\,\,|\,\, R_1,\dots, R_m,R_{m+1},\dots, R_{m+n}
\rangle
\]
This  and (\ref{relation1}) implies that there is a canonical surjection
$\phi:\pi_1(\BC_L^2-D_0)\to \pi_1(\BC_L^2-C_\tau)$.\qed
\begin{Remark}
Though we are mainly concerned in the case $j\ge 2$, the assertion for
 $j=0$
gives another proof of  Theorem 3.
\end{Remark}
\begin{Corollary}Under the same assumption,
we have the divisibility
$\De_{C_\tau}(t;L)\,|\, \De_{D_0}(t;L)$.
\end{Corollary}
Taking  a generic line $L_\infty$ for $C_\tau$, $\tau$ small and
$D_0\cup L$, we apply Theorem \ref{Tangential} to obtain:
\begin{Corollary}Under the same assumption,
\[
 \sqrt{\De_{C_\tau}(t)}\,\,|\,{\De_{D_0\cup L^j}(t)}
\]
\end{Corollary}
{\em Proof.}
First note that $\De_{C_\tau}(t)\,|\, \De_{C_\tau}(t; L)$
and $\De_{C_\tau}(t;L)\,|\, \De_{D_0}(t;L)$ by Corollary 15.
Secondly
$\De_{D_0\cup L_\infty^{-d}}(t;L)=\De_{D_0\cup L^j}(t)$ by
 Theorem \ref{Tangential}.
Lastly  we have
\[
 \sqrt{\De_{D_0\cup L_\infty^{-d}}(t;L)}=
\gcd({\De_{D_0}(t;L)\times(t-1)},t^d-1)
\]
The conclusion  is now immediate from these observation,
as the factor $(t-1)$ does not appear in $\De_{C_\tau}(t)$
by the irreducibility \cite{OkaSurvey}.
\qed
\subsubsection{Examples of line degenerations of order 1.}
In \cite{Reduced},
we have classified configurations of reduced sextics of torus type.
Among them, there are sextics of torus type with components $B_5+L$
where $B_5$ is an irreducible quintic. In fact, each of them is a
 degeneration of irreducible sextics of torus type.  We give one such
 example.
$D_0=B_5$ has one $A_{11}$ singularity at $O=(0,0)$ and $L=\{Z=0\}$
and it is a flex tangent.
A generic irreducible sextic $C_t$ has $[A_{11}+A_5]$ as singularities.
By Theorem \ref{surjectivity}, we know that 
$\De_{C_t}(t;L)\,|\,\De_{D_0}(t;L)$. As $\De_{C_t}(t;L)$ is divisible by
 the generic Alexander polynomial, which is $t^2-t+1$ by
 \cite{OkaAtlas},
we conclude $t^2-t+1\,|\,\De_{D_0}(t;L)$.

\begin{eqnarray*}
&D_0:\quad {x}^{4}y-8\,{x}^{2}{y}^{2}+7\,{y}^{3}-4\,{x}^{3}y-2\,{x}^{2}y+{y}^{2}+
{x}^{4}+2\,{x}^{5}\\
&\qquad\qquad +9\,{y}^{4}-6\,{x}^{3}{y}^{2}+2\,{y}^{2}x+6\,{y}^{3}
x=0\\
&C_s:\,\,\left( -{x}^{2}+y \right) ^{3}+ \left( -{x}^{2}-{x}^{3}+y+3\,{y}^{2}+
xy+s\,{x}^{2}y \right) ^{2}=0
\end{eqnarray*}

\subsubsection{ Line degenerations of order 2.}
We consider quartics which are non-reduced line degenerations of sextics.
As a quartic $D$, we can take quartics with configuration 
(a) $\Si(D)=2A_2$, or (b)  $\Si(D)=2A_2+A_1,\, 3A_2$
 (with an outer singularity
$A_1$ or $A_2$)
or (c) $\Si(D)=A_5,\,E_6$ or $\Si(D)=A_2+A_3,\,A_6$ (with a wild
inner singularity).

We only explain here the case $(b)$ with $\Si(D)=3A_2$. 
The other case will be explained
systematically
in \S 5.4. 
Three cuspidal quartic is very special. 
As the degenerated
line, we can take 
any one of  a simple tangent line
or a tangent cone at a cusp or a unique bi-tangent line.
This is not the case for other quartics listed above.
 %
 The following family gives line degeneration of order 2 of sextics of
torus type defined  by
$f_2^3+g_3^2=0$ into a quartic with 3 $A_2$ and a line $\{Z=0\}$
with multiplicity 2.

\vspace{.2cm}\noindent
(a-1) $L$ is a {\bf simple tangent line}. 
$\{C_s,\,s\in \De\}$ is a family with   an outer $A_2$ singularity at $Q=(0,1,0)$ and $L=\{Z=0\}$
is a simple tangent line of the quartic $D_0$. This degeneration is not a
line-degeneration of
torus curve which we study in the next section.
\begin{eqnarray*}
&D_0:\, {y}^{4}-2\,{y}^{2}+3\,x{y}^{2}-3/4\,{x}^{2}{y}^{2}+1-3\,x+3\,{x}^{2}-{
x}^{3}=0\\
&C_s^{(1)}:\,\left( y-3/2\,xy-{y}^{3}+\,s\,{x}^{3} \right) ^{2}+
 \left( -{y}^{2}+1-x
 \right) ^{3}=0
\end{eqnarray*}

\begin{figure}[htb]
\setlength{\unitlength}{1bp}
\begin{picture}(600,150)(-100,0)
 \put(100,160){\special{epsfile=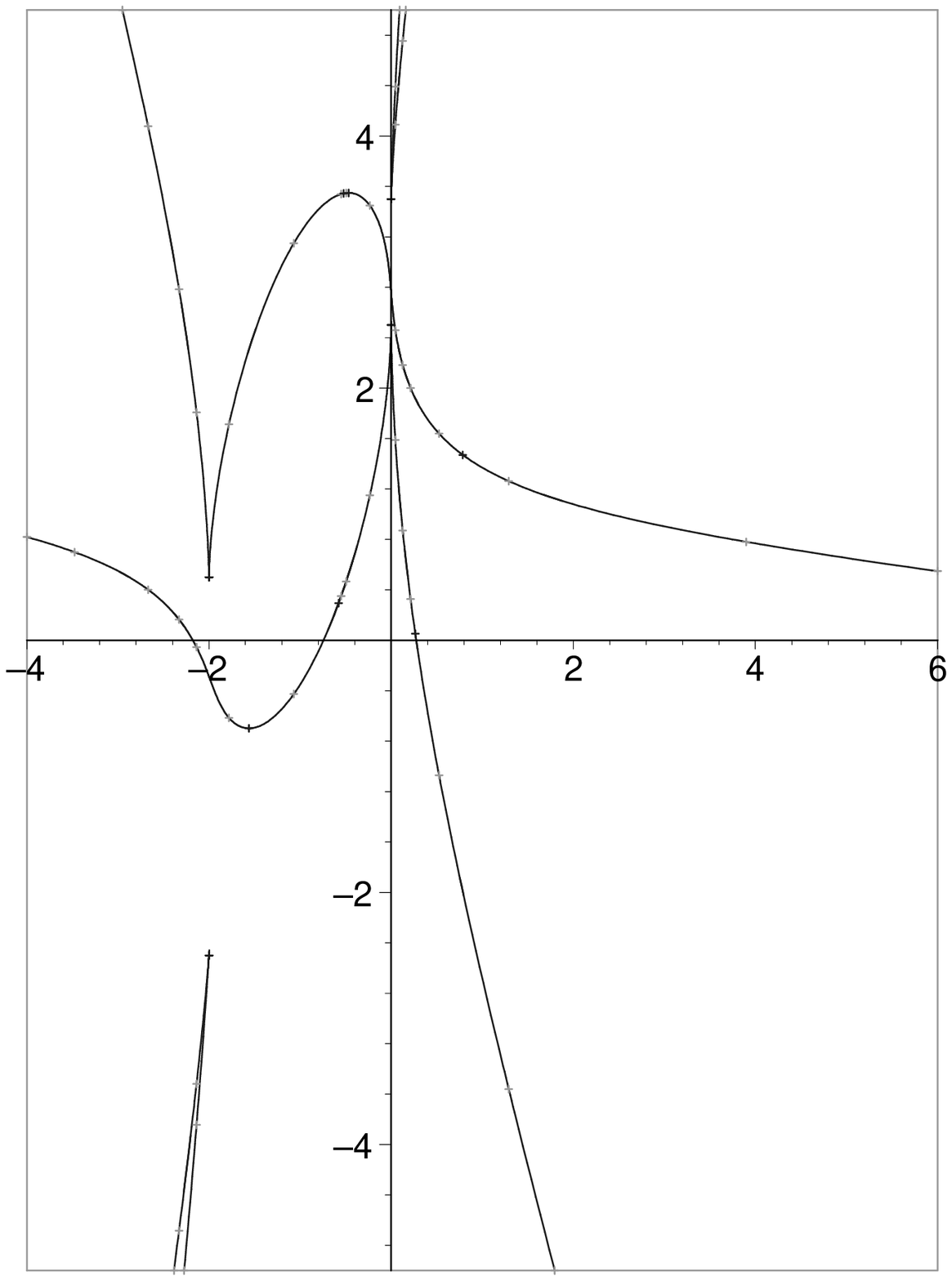 vscale=0.5
 hscale=0.6}}
\end{picture}
\vspace{3cm}
\caption{$C_s^{(3)},\,s=-1/3,\,$sextics with $8A_2$}
\label{Quartic-3A2}
\end{figure}

\noindent
(a-2)  $L$ is a {\bf  tangent cone of an $A_2$}. 
$\{C_s^{(2)},\,s\in \De\}$ is a family with  an $A_1$ singularity at $Q=(0,1,0)$ and $L=\{Z=0\}$
is a  tangent cone of an $A_2$ singularity.
\begin{eqnarray*}
&D_0:\,-3\,x{y}^{3}+{y}^{4}+3\,xy+1+3\,{x}^{2}{y}^{2}-5/4\,{x}^{3}y+3/16\,{x}
^{4}-3/4\,{x}^{2}-2\,{y}^{2}=0\\
&f_2(x,y,s)=1+ \left( -1/4+s \right) {x}^{2}+xy-{y}^{2} \\
&g_3(x,y,s)=-{y}^{3}+3/2\,x{y}^{2}+ \left( -3/4\,{x}^{2}+1 \right) y+1/8\,{x}^{3}
\end{eqnarray*}
\noindent
(a-3)  $L$ is the {\bf  bi-tangent line}. 
As is observed in Proposition \ref{Exceptional},
the Alexander polynomial has multiplicity 2 for $(t^2-t+1)$.
Thus this case is exceptional for the tangential
Alexander polynomial of quartic with $3A_2$. To expalin
this we consider
 a family of sextics of torus type $\{C_s^{(3)},\,s\in \De\}$  with $8A_2$,
where two $A_2$ are outer singularities and  they are located at 
at $R=(0,5/2,1)$ and $R'=(-2,1/2,1)$. An inner $A_2$ is at
$Q=(0,7/2,1)$.
 The line $L=\{X=0\}$
is a  bi-tangent line and the limit line degeneration. 
In an affine equation, we can define them as
 \begin{multline*}
D_0:\,
x{y}^{3}+{\frac {1}{256}}\, \left( 19\,{x}^{2}-12\,{y}^{2}+75\,x+75+6
\,xy \right) ^{2},\,C_s^{(3)}=\{f_2^3+g_3^2=0\}\\
f_2(x,y,s)=
\frac 1{4(s-1)}\,(20\,ysx+8\,{s}^{2}xy-8\,{s}^{3}xy-70\,s+4\,{s}^{
3}{x}^{2}-20\,{s}^{3}y-21\,{s}^{2}+21\,{s}^{3}-4\,xy\\
-4\,{s}^{2}{x}^{2}+20\,{s}^{2}y
-8\,{y}^{2}s-4\,{s}^{2}{y}^{2}+4\,{s}^{3}{y}^{2}+48\,ys-8
\,s{x}^{2}-48\,sx+18\,{s}^{3}x-18\,{s}^{2}x)\\
g_3(x,y,s)
=\frac 1{16(s-1)^2)} (552\,ysx+36\,x{y}^{2}{s}^{5}-90\,y{s}^{2}{x}^{2
}+138\,ys{x}^{2}-42\,y{s}^{3}{x}^{2}+36\,{y}^{2}{s}^{3}x-84\,{y}^{2}x{
s}^{4}\\
+36\,{y}^{2}{s}^{2}x-108\,{y}^{2}sx-336\,{s}^{2}xy-96\,{s}^{3}xy
-75\,x+516\,x{s}^{4}y+120\,{x}^{2}{s}^{4}y-228\,{s}^{5}xy-48\,{s}^{5}{
x}^{2}y\\
-525\,s+24\,{y}^{3}s+333\,{s}^{5}x-75\,{x}^{2}+117\,{s}^{3}{x}^
{2}+18\,{s}^{3}y+420\,{s}^{2}-161\,{s}^{3}+21\,{s}^{5}{x}^{3}-741\,x{s
}^{4}\\
+12\,x{y}^{2}-13\,s{x}^{3}-518\,{s}^{4}+259\,{s}^{5}+16\,{s}^{4}{
y}^{3}-270\,{s}^{5}y-354\,{s}^{4}{x}^{2}+177\,{s}^{2}{x}^{2}-288\,{s}^
{2}y\\
-8\,{s}^{5}{y}^{3}-8\,{s}^{3}{y}^{3}-168\,{s}^{4}{y}^{2}-6\,{x}^{2
}y+540\,{s}^{4}y-19\,{x}^{3}+8\,{s}^{2}{x}^{3}-204\,{y}^{2}s+48\,{s}^{
2}{y}^{2}\\
+36\,{s}^{3}{y}^{2}+144\,{s}^{5}{x}^{2}+570\,ys-213\,s{x}^{2}
-645\,sx+15\,{s}^{3}x+543\,{s}^{2}x+36\,{s}^{3}{x}^{3}-57\,{s}^{4}{x}^
{3}
+84\,{y}^{2}{s}^{5})
\end{multline*}
As the Alexander polynomial of sextics of torus type with $8A_2$
is given $(t^2-t+1)^2$ (\cite{OkaAtlas}), this explains 
$\De_{D_0}(t;L)=(t^2-t+1)^2$.
Four inner $A_2$'s
are not visible in the Figure \ref{Quartic-3A2}.
It is quite interesting to study how the family degenerates into
$D_0+2L$. Observe that two cusps of $D_0$ are not  real points
and $L=\{x=0\}$ is the bitangent line of $D_0$.
Note that (a-1) and (a-2) can not be a line-degenerated torus curve
in the sense of the next subsection.

\vspace{.5cm}
\noindent

\subsection{Line degeneration of curves of torus type }
We consider a pair of coprime  positive integers $p,\,q>1$ and 
consider the curves of $(p,q)$-torus curve:
\[
 C_{p,q}:\quad f_p(X,Y,Z)^q+g_q(X,Y,Z)^p=0
\]
where $f_p,g_q$ are polynomials of degree $p,\,q$ respectively.
Consider the special case that
\[
 f_p(X,Y,Z)=f_{p-a}(X,Y,Z)\times Z^a,\quad
g_q(X,Y,Z)=g_{q-b}(X,Y,Z)\times Z^b
\]
where $f_{p-a},\,g_{q-b}$ are homogeneous polynomials of degree $p-a,\, q-b$ respectively and
$0<a<p$ and $0<b<q$. Assume  for example that
$ aq\,\ge\, bp$
and factoring $Z^{bp}$ from $f$, we have a curve 
\begin{eqnarray}\label{deg1}
D:\quad g(X,Y,Z)\,=\,f_{p-a}(X,Y,Z)^q\,Z^{qa-pb}\,+\,g_{q-b}(X,Y,Z)^p=0
\end{eqnarray}
We call  a curve $D$ {\em a  line-degenerated torus curve of type $(p,q)$}
and we call the line $L=\{Z=0\}$  the {\em   limit line of
the  degeneration}.
Note that the degree of $D$  is $pq-bp$.
The simplest case is $a=b=1$ and
\[
 D:\quad  g(X,Y,Z)\,=\,f_{p-1}(X,Y,Z)^q Z^{q-p}+ g_{q-1}(X,Y,Z)^p=0,\quad p<q
\]
\begin{Theorem}\label{torus-degeneration}
For a degenerated torus curve $D$ of type $(p,q)$ defined  (\ref{deg1}),
there is a family of line degeneration $f(X,Y,Z,t)$ of order $bp$ 
such that $ f_0(X,Y,Z)=g(X,Y,Z)\, Z^{bp}$ and 
each curve $C_t: \,f(X,Y,Z,t)=0$
passes through a fixed point $Q\in L$ and the multiplicity
of $C_t$ at $Q$ is  $pb$ for each $t\ne 0$.
\end{Theorem}
{\em Proof.} We may assume that $Q=(0,1,0)$ is not on $D$.
Let $h(X,Y)$ be a homogeneous polynomial of degree $q-b$ with 
$h(0,Y)\ne 0$.
We put 
\[
 f_p(X,Y,Z)=f_{p-a}(X,Y,Z)Z^a,\quad
g_q(X,Y,Z,t)=g_{q-b}(X,Y,Z)Z^b\,+\, t\,X^b h(X,Y)
\]
and put  $$C_t:\, f(X,Y,Z,t):=f_p(X,Y,Z)^q+g_q(X,Y,Z,t)^p=0
$$
We can easily see that $C_t:f(X,Y,Z,t)=0$ passes through $Q$
and the multiplicity of $(C_t,Q)$ is $pb$, as the local equation at $Q$
is given by
\[
 f(x',1,{z'},t)=f_{p-a}(x',1,{z'})^q {z'}^{qa}\,+\, \left(g_{q-b}(x',1,{z'}){z'}^b\,+\, t\,{x'}^b h(x',1)\right)^p=0
\]
where $x'=X/Y,\,{z'}=Z/Y$.
The affine equation of $C_t$ in $\BC_L^2$ is given by
\[
 f(x,y,1,t)=f_{p-a}(x,y,1)^q\,+\, \left(g_{q-b}(x,y,1)\,+\, t\,x^b h(x,y)\right)^p=0.
\]
where  $x=X/Z,\,y=Y/Z$.
We  see that
$\degree_y\, f(x,y,1,t)=pq-pb$.\qed

Put $s=\gcd(p,q)$.
As  the generic Alexander polynomial of $(p,q)$-torus curve
of degree $pq$ is divisible by
$\De_{p,q}(t):={(t^{pq/s}-1)^s(t-1)}/((t^p-1)(t^q-1))$, we get
\begin{Corollary}\label{divisibility-torus} Let $D$ be as above and let
 $L=\{Z=0\}$ be the
limit
 line
of the degeneration.  Then
$\De_D(t;L)$ is divisible by $\De_{p,q}(t)$.
\end{Corollary}
\subsubsection{Singularities of line-degenerated torus curves}
We consider the curve defined by (\ref{deg1})
\begin{eqnarray*}
 D:\quad g(X,Y,Z)\,=\,f_{p-a}(X,Y,Z)^q\,Z^{qa-pb}\,+\,
g_{q-b}(X,Y,Z)^p=0
\end{eqnarray*}
Suppose that $P\in D$ is a singular point.  $P$ is called 
{\em an inner singularity} (respectively {\em outer}) if 
$f_{p-a}(P)=g_{q-b}(P)=0$ (resp. $g_{q-b}(P)\ne 0$). $P$ is called  {\em wild} if $P$ is also
on the limit line $Z=0$.
The following describes the type of  inner non-wild singularities.
\begin{Lemma} (\cite{BenoitTu})\label{sing-type1}
Let $C$ be a curve of torus type
\[
 C:\quad f(x,y)^p+g(x,y)^q=0,\, p<q
\]
and assume that $f(0,0)=g(0,0)=0$ and 
the curves $f(x,y)=0$ is smooth at $O$. Let $\nu$
be the local intersection number of $f(x,y)=g(x,y)=0$ at $O$.
Then the singularity $(C,O)$ is topologically isomorphic to the
 Brieskorn singularity
\[
 B_{p,q\nu}:\,\, y^p\,+\, x^{q\nu}=0
\]
\end{Lemma}
The singularity is more complicated when $f(x,y)=0$
is singular at $O$. The description of wild singularities is also
more complicated in general.
\subsection{Examples of line degeneration of  torus curves}
\subsubsection{Cubic}
A cuspidal cubic $Q$ can be understood as a line-degenerated (2,3)-torus
curve  of order 3  by taking 
\begin{eqnarray*}
&f_2(x,y)=yx,\, g_3(x,y)=g_1(x,y)\, x^2\\
& Q:\,\, y^3\,+\, g_1(x,y)^2x=~0.
\end{eqnarray*}
The limit line of the degeneration is $L=\{x=0\}$. This explains that $\De_Q(t;L)=t^2-t+1$.

\subsubsection{Quartic}
We give further quartics which can be a line-degenerated
$(2,3)$-torus curve.
We consider the quartics of the form:
\begin{eqnarray}\label{gen-quartic}
C_0:& \,\, (f_2(x,y)\,x)^2\,+\, (y\,x)^3=0\\
D:&\,\, g(x,y))=f_2(x,y)^2\,+\, y^3\, x=0
\end{eqnarray}
where $f_2(x,y)$ is a polynomial of degree
2 and  the limit line  of degeneration is chosen to be $\{x=0\}$.
In general, $D$ has two inner $A_2$ singularities
at $y=f_2(x,y)=0$. If $y=0$ is tangent to  the conic $C_2:=\{f_2(x,y)=0\}$,
the singularity is an $A_5$. If moreover $C_2$ degenerates  into
two lines, the singularity is  an $E_6$-singularity.
The limit line $L$ is a bi-tangent line at 
$\{x=0\}\cap C_2$.  If $L$ is tangent to the conic $C_2$,
$D$ obtains a flex of order 2 and  $L$ is the flex tangent line.
Further more we can put one outer singularity, either $A_1$ or $A_2$.
There are two more configurations  which can be a line degeneration of
torus curves: $A_2+A_3$ and $A_6$.
For these singularities, we have to consider wild inner singularities.
We have already studied most of these quartic and their Alexander invariants
in \S \ref{example-quartic}.  Theorem \ref{surjectivity} and 
Theorem \ref{torus-degeneration} explains our previous computations.
\subsubsection{ A quartic with $2A_2$}
For a generic  quartic $D$ with $\Si(D)=2A_2$, its dual curve $\check D$ is a sextic
with 
$8A_2+A_1$,  in particular, $D$ has one bi-tangent line.
\begin{eqnarray*}
D:\,\,g(x,y)= \left( {y}^{2}-1+{x}^{2} \right) ^{2}+{y}^{3}x=0
\end{eqnarray*}
The bi-tangent line can degenerate into a flex tangent line of order 2
so that $\check D$ has $6A_2+E_6$.
\[
 D:\,\,g(x,y)= \left( {y}^{2}-2\,y-2\,{x}^{2}+x+1 \right) ^{2}+{y}^{3}x=0
\]
See Figure \ref{Quartic-2A2}, \S 5.9 for graphs of these quartics.
\subsubsection{A quartic with $2A_2+A_1$}
A generic  quartic $D$  with $\Si(D)=2A_2+A_1$
 has two flexes and one bi-tangent
line (i.e., the configuration space is self dual). Degenerated quartic
$D'$ has one flexes of order 2 and no bi-tangent line as
we have seen before.
\begin{eqnarray*}
D:\,\, \left( {y}^{2}+ \left( \frac 92\,x-1 \right) y+\frac 32\,{x}^{2}+x-1 \right) ^{
2}+{y}^{3}x=0
\end{eqnarray*}
\begin{eqnarray*}
D':\,\,
\left( {y}^{2}+ \left( {\frac {17}{6}}\,x-\frac 73 \right) y+{\frac {23}{6
}}\,{x}^{2}-\frac {13}3 \,x+1 \right) ^{2}+{y}^{3}x=0
\end{eqnarray*}
See Figure \ref{Quartic-2A2A1}, \S 5.9 for graphs of these quartics.
\subsubsection{A quartic with $A_5$}
A generic  quartic $D$  with $\Si(D)=\{A_5\}$ has 6 flexes and one bi-tangent
line, thus the dual curve is a sextic with
$6A_2+A_5+A_1$. 
 Degenerated quartic
$D'$ has 4 flexes  and one flex of order 2 and thus the dual curve is a
sextic with $4A_2+E_6+A_5$. Tokunaga has studied a certain dihedral
covers branched along these quartics \cite{Tokunaga-Kyoto}
\begin{eqnarray*}
D:\,\, \left( {y}^{2}-yx-{x}^{2}+2\,x-1 \right) ^{2}+{y}^{3}x=0\\
D':\,\,\left( {y}^{2}-2\,y+{x}^{2}+2\,x+1 \right) ^{2}+{y}^{3}x=0
\end{eqnarray*}
See Figure \ref{Quartic-A5}, \S 5.9 for graphs of these quartics.
\subsubsection{A quartic with $E_6$}
A generic  quartic $D$  with $\Si(D)=\{E_6\}$ has two flexes and one bi-tangent
line. Degenerated one $D'$ has one flexes of order 2 and no bi-tangent line as
we have seen before.
\begin{eqnarray}
&D:\,\,\left( {y}^{2}-1-{x}^{2}+2\,x \right) ^{2}+{y}^{3}x=0
& D':\,\, \left(x-y+1\right) ^{4}+
{y}^{3}x=0
\end{eqnarray}
Note that the last quartic $D'$ with $E_6$ can be also considered as a line
degeneration of (3,4)-torus curves as
\[
 \left((x-y+1)^4+y^3 x\right) x^8=
((x-y+1)x^2)^4+(y x^3)^3
\]
This explains that $\De_{D'}(t;L)=(t^2-t+1)(t^4-t^2+1)$ with
$L=\{x=0\}$.
See Figure \ref{Quartic-E6}, \S 5.9 for graphs of these quartics.
\subsubsection{quartics with $A_2+A_3$ and $A_6$}
We start the general form
\[
 C:\,\,f_2(x,y)^2\,+\, y^3\,x=0
\]
We assume that $O=(0,0)$ is a wild inner singularity. Thus 
$f_2(0,0)=0$.

The configuration $A_2+A_3$ is obtained when 
$f_2(0,0)=0$ and $f_2(x,y)=0$ is not tangent to $y=0$
as
$(C,O)\cong A_3$. The limit line is tangent at one smooth point and also
passing at $O$. For example,
\[
C:\,  \left( -{y}^{2}+y-{x}^{2}+x \right) ^{2}+{y}^{3}x=0
\]
If we make $f_2(x,y)=0$ is tangent to $y=0$ at $O$, $(C,O)\cong A_6$.
An example is given as follows.
\[
 C:\, \left( {y}^{2}+y+{x}^{2} \right) ^{2}+{y}^{3}x=0
\]
See Figure \ref{Quartic-A3A2}, \S 5.9 for graphs of these quartics.
\subsubsection{Quintics as line-degenerations}
We consider  (2,5)-torus curve $C: f(x,y)=f_5(x,y)^2+g_2(x,y)^5=0$ of degree 10 which is degenerated as 
\[
 f_5(x,y)=f_2(x,y) x^3,\, g_2(x,y)=yx
\]
Then we get a quintic
\begin{eqnarray}\label{quintic-general}
 D:\,\,f_2(x,y)^2\,x+ y^5=0
\end{eqnarray}
In general, $D$ has 2 $A_4$ singularities at $f_2(x,y)=y=0$ and it has a
flex of
order 3 at $O$ with the tangent line $L=\{x=0\}$.
As a special case where the conic $f_2(x,y)=0$ is tangent to $y=0$, we
get
one $A_9$ singularity:
\[
 D: \,\,
 \left( {y}^{2}+xy+{x}^{2}-2\,x+1 \right) ^{2}x+{y}^{5}=0
\]
If $f_2(x,y)=0$  is  two lines intersecting  on $y=0$,
the singularity is locally topologically isomorphic to  $C_{5,5}$ 
in the notation \cite{Pho}:
\[
 C_{5,5}:\, y^5+x^2y^2+x^5=0
\]
 If $f_2=0$ is a line with multiplicity 2,
the singularity becomes  $B_{4,5}$ 
singularity which is locally defined as $\,y^5+x^4=0$.

There are two other possibilities. A quintic as a line degeneration of
torus curves of type (3,5): take $g_3=yx^2,\,f_5=f_1(x,y)x^4$.
Then we get a quintic
\[
 Q:\quad y^5+x^2\,f_1(x,y)^3=0
\]
The quintic $Q$ has one $E_8$ singularity and a $A_4$ singularity on the
limit line $x=0$. The limit line is also the tangent line of the singularity $A_4$.

Another possibility is as 
a line degeneration of
torus curves of type (4,5): take $g_4=yx^3,\,f_5=g_1(x,y)x^4$.
Then we get a quintic
\[
 Q':\quad y^5+x\,g_1(x,y)^4=0
\]
$Q'$ has one $B_{5,4}$-singularity and the limit line is $x=0$. This
quintic can be considered as a degeneration of (\ref{quintic-general})
when
$f_2(x,y)\to g_1(x,y)^2$.
\subsection{Sextics as line degenerations}
Sextics as line degenerations can be either
 from (2,5)-torus curves, or  from (3,5)-torus curves
or   from (3,4)-torus curves.
The sextics from (2,5)-torus curves take the form:
\begin{eqnarray}
C:\,g(x,y)=f_3(x,y)^2+y^5\,x=0.
\end{eqnarray}
Generically $C$ has $3$ $A_4$ singularities and the  degeneration line
$L=\{x=0\}$ is a tri-tangent line.
By the degeneration of the intersection $y=f_3(x,y)=0$, we may have also 
either $A_4+A_9$ or $A_{14}$. If the cubic $f_3(x,y)=0$ has  a node or
cusp, the singularity 
 becomes more complicated.

Sextics from (3,5) torus curves take the form:
\begin{eqnarray}
C':\, f_2(x,y)^3+y^5\,x=0.
\end{eqnarray}
Generically $C'$ has 2 $B_{3,5}$ singularities and degeneration line is
 a bi-tangent line at two flex points. If $y=0$ is tangent to the conic
 $f_2(x,y)=0$,
the singularity is $B_{3,10}$.

The curves $C:\,f_2(x,y)^3+ y^4\,x^2=0$
can be considered as a line degeneration of (3,4)-torus curves
but at the same time, it is a torus curve of type (2,3). Generically $C$
has $2E_6+2A_2$.
Thus by Corollary \ref{divisibility-torus}, the Alexander polynomial
$\De_C(t;L)$
is divisible by $(t^2-t+1)(t^4-t^2+1)$.
Note that $\De_C(t)=t^2-t+1$ by \cite{Oka-Pho1}.

\subsection{Flex degenerations}
Let us consider  $C_\tau$ is a family of irreducible curves
in the configuration space $\cM(\Si;d)$ with two marked
flex points $P_\tau,\,Q_\tau\in C_\tau$ of order 1 for $\tau\ne 0$ and assume that 

(a) $P_\tau,\,Q_\tau\to  P_0$ when $ \tau\to 0$ and $P_0$ is a flex point of order 
$2$
of  $C_0\in \cM(\Si;d)$.

(b) The intersection $ T_{P_\tau} C_\tau\cap  T_{Q_\tau} C_\tau\cap C_\tau$
 is empty for $\tau\ne 0$.

\begin{Theorem}\label{divisibility3} Consider the 
degeneration:
$C_\tau+L_{1,\tau} +  L_{2,\tau}  \to C_0+2L$. Then
we have the divisibility of Alexander polynomials:
\[
 \De_{C_\tau\cup L_{1,\tau}\cup L_{2,\tau}}(t)\,|\,\De_{C_0\cup L^2}(t)\times (t-1).
\]
Here $L_{1,\tau}=T_{P_\tau} C_\tau$,  $ L_{2,\tau}=T_{Q_\tau} C_\tau$
 and $L=T_{P_0} C_0$.
\end{Theorem}

{\em Proof.} 
First we may assume that $C_\tau=\{f(x,y,\tau)=0\}$,
 $P_0=(0,0)$ and the tangent line of $C_0$ at $P_0$ is defined by
$y=0$. Take a  generic line $L_\infty$ and we work in $\BC^2=\BC_{L_\infty}^2$.
Taking a presenting polydisk $\De_{\al,\be}$ for
$C_0\cup L$ with pencil line $L_\eta,\,\eta\in \BC$, let $\eta_1,\dots, \eta_m$ be parameters corresponding to
the singular pencil lines with $\eta_1=0$.
 Fix a small $\eps>0$ to see the monodromy
relations
along
$|\eta-\eta_i|=\eps$. We take generators in a fixed generic line
$L_{\eta_0}, \, |\eta_0|=\eps$ (see Appendix).
Then we get a  presentation:
\begin{eqnarray}\label{rep-L}
\pi_1(\BC^2-C_0\cup L)=\langle g_1,\dots, g_d,h\,|\, R_1,\dots, R_k\rangle
\end{eqnarray}
We consider the Alexander matrix $M_0$ with respect to the weight
function $\theta_2$  which  has  weight 2  for
 $L$.
Take a positive number $\de$ so that 
\[
 C_\tau\pitchfork L_\eta,
\quad \forall \tau,\,|\tau|\le \de,\,\forall\eta,\,|\eta-\eta_i|=\eps,\, 
i=1,\dots,m
\]
We may also assume that $\De_{\al,\be}$ is admissible for 
$C_\tau,\, |\tau|\le \de$.
Next, we consider $C_\tau\cup L_{1,\tau}\cup L_{2,\tau}$ for sufficiently
small $\tau$ in the above sense.
 In the generic fiber $L_\eta$, the intersection
$L_\eta\cap L$ are  bifurcated into two points $L_\eta\cap L_{1,\tau}$
and $L_\eta\cap  L_{2,\tau}$, but they  are observed  only  with a
microscope and they move exactly as a twin satellite along 
$|\eta-\eta_i|=\eps$. For the presentation of $\pi_1(\BC^2-C_\tau\cup
L_{2,\tau}\cup L_{2,\tau})$,
we need two generators $h_1,\,h_2$ presented by lassos 
for the lines $L_{2,\tau}$ and
$L_{2,\tau}$ instead of one $h$. However we can understand as
$h=h_1h_2$. For the further detail about the choice of generators, 
see Appendix.
This implies that there are canonical  homomorphisms $\psi,\,\Psi$
which make  the next diagram commutative.
\[
 \begin{matrix}
F(d+1)&\mapright{}&\pi_1(\BC^2-C_\tau\cup
  L)&\mapright{\theta_2}&\BZ\\
\mapdown{\Psi}&&\mapdown{\psi}&&\mapdown{\id}\\
F(d+2)&\mapright{}&\pi_1(\BC^2-C_\tau\cup L_{1,\tau}\cup L_{2,\tau})&
\mapright{\theta_{sum}}&\BZ
\end{matrix}
\]
Here $\Psi$ is defined on generators as $g_i\mapsto g_i,\,h\mapsto
h_1h_2$
and $\psi$ is canonically induced by $\Psi$.
The monodromy relations
$R_1,\dots, R_k$ remains the same along $|\eta-\eta_j|=\eps$.
 This  means that
the relation 
$R_j'$ remains true where $R_j'$ is obtained simply
substituting the letter $h$ by $h_1h_2$.
To get a complete relations, we have to add some more relations, say
$S_1,\dots, S_\ell$ along
the singular pencil lines. Among them, we can assume that
$$S_1=g_dh_2g_d\inv h_2\inv
$$
 which is the relation at the transverse intersection 
$L_{2,\tau}\cap C_{\tau}$ near $P_0$. (In the appendix, we will explain
this situation more.)
Thus the presentation is given as
\begin{eqnarray}\label{rep-L1L2}
\pi_1(\BC^2-C_\tau\cup L_{1,\tau}\cup L_{2,\tau})=
\langle g_1,\dots, g_d,h_1,h_2\,|\, R_1',\dots, R_k',\, S_1,\dots, S_\ell\rangle
\end{eqnarray}
Let $G_\infty$ be the subgroup of $G_\tau:=\pi_1(\BC^2-C_\tau\cup L_{1,\tau}\cup
L_{2,\tau})$
generated by $g_1,\dots, g_d$ and the product $ h_1h_2$. Then $\psi$ is a surjection on
$G_\infty\subset G_\tau$.
For $C_\tau\cup L_{1,\tau}\cup L_{2,\tau}$, we consider the summation
homomorphism $\theta_{sum}$. Let $\gamma_2,\,\gamma_\tau$ be the
 ring homomorphisms corresponding to $\theta_2,\, \theta_{sum}$:
\[
 \gamma_2:\,\BC(F(d+1))\to \BC[t,t\inv],\,\, \gamma_\tau:\BC(F(d+2))\to \BC[t,t\inv]
\]
Note that the following diagrams are commutative.
\[
 \begin{matrix}
\ga_2:\BC(F(d+1))&\mapright{}&\BC(\pi_1(\BC^2-C_\tau\cup
  L))&\mapright{{\theta_2}_*}&\BC(\BZ)=\La\\
\mapdown{\Psi_*}&&\mapdown{\psi_*}&&\mapdown{\id}\\
\ga_\tau:\BC(F(d+2))&\mapright{}&\BC(\pi_1(\BC^2-C_\tau\cup L_{1,\tau}\cup L_{2,\tau}))&
\mapright{{\theta_{sum}}_*}&\BC(\BZ)=\La
\end{matrix}
\]
Now we consider the Alexander matrix $M_\tau$ of  $C_\tau\cup
L_{2,\tau}\cup L_{2,\tau}$.
We consider the row corresponding to 
the relation $R_i'$. It is a word of $g_1,\dots, g_d$ and $h_1h_2$.
By the definition of $\theta_2$, we can see easily that 
\[
 \gamma_2(\frac{\partial R_i}{\partial g_j} )=\gamma_\tau(\frac{\partial R_i'}{\partial g_j} )
\]
The $(d+1)$-th column split into two columns,
which correspond to the Fox differentials
$\frac{\partial}{\partial h_i},\,i=1,2$.
As $h=h_1h_2$, we obtain
\begin{eqnarray}
(\gamma_\tau(\frac{\partial R_i'}{\partial  h_1}),\gamma_\tau(\frac{\partial R_i'}{\partial h_2}))=
( \gamma_2(\frac{\partial R_i}{\partial h} ), t \times \gamma_2(\frac{\partial R_i}{\partial h})
\end{eqnarray}
Note also that 
\[
 (\gamma_\tau(\frac{\partial S_1}{\partial
 h_1}),\gamma_\tau(\frac{\partial S_1}{\partial h_2}))=
(0,t-1)
\]
Let $M_\tau'$ be the matrix obtained by adding
 $(-t)\times (d+1)$-th column to $(d+2)$-th 
column so that the last column is zero up to $k$-th row.
$M_\tau'$ is written as
\[
 \left(
\begin{matrix}
M_0&\vec 0\\
\bfw&(t-1)\\
N_1&\vec v
\end{matrix}
\right)
\]
where $M_0$ is the Alexander matrix of $C_0\cup L$
and $\bfw$ comes from the differential of $S_1$,
and the other terms
 $N_1,\, \vec v$  are coming from 
$S_j,\,j\ge 2$. Thus for any $(d-1)\times (d-1)$-minor $A$ of
 $M_0$,
we associate  $d\times d$-minor $A'$ of $M_\tau'$, by adding $(k+1)$-row and the last
column.
Then  the corresponding determinant is equal to $\det(A)\times (t-1)$.
Thus the assertion  follows from the Fox calculus definition of the
Alexander polynomial.\qed

\subsection{Appendix} In this appendix, we will explain the existence
of the relation $S_1$ in the proof of Theorem \ref{divisibility3}.
First we may assume that $C_\tau=\{f(x,y,\tau)=0\}$,
 $P_0=(0,0)$ and the tangent line of $C_0$ at $P_0$ is defined by
$y=0$. Changing the scale and using Implicit function theorem, we
may assume that $C_0$ is defined by
$y=\phi_0(x)$ where $\phi_0(x)=x^4+\text{higher terms}$
in the polydisk $\De_{1,1}=\{(x,y);|x|,\,|y|\le 1\}$.
This follows from the assumption that $O$ is a flex of order 2
of $C_0$. We consider the
pencil
lines $x=\eta$.
Now we consider
$C_\tau$. Assume that
$C_\tau$ is defined by 
$y=\phi_\tau(x)$ in $\De_{1,1}$.
 Write $\phi_\tau(x)=\sum_{\nu}c_\nu(\tau)x^\nu$.
First observe that $|c_4(\tau)-1|\ll 1$ by continuity.
 In the
parametrization, the flex points are defined by 
$\{(\alpha,\beta)| \phi_\tau''(\alpha)=0\}$.
Thus by Rouch\'e's principle, we see that there is two flex points
which bifurcate from $P_0$. They corresponds to the 
roots of $\phi_\tau''(x)=0$ in $|x|\le 1$, say $x=\alpha_1(\tau),
\alpha_2(\tau)$. Thus $L_{1,\tau},\,L_{2,\tau}$ corresponds to the
tangent line at these flex points. By B\'ezout's theorem and the local
stability of intersection numbers, there is one transverse intersection
point of
$L_{i,\tau}\cap C_\tau$ and we put them 
 $Q_i=(\beta_i,\phi_\tau(\beta_i))$ for $i=1,2$.
Note that $\beta_i\to 0$  when $\tau\to 0$.
Thus the local singular pencils for 
$\widetilde C_\tau=C_\tau\cup L_{1,\tau}\cup L_{2,\tau}$ is bifurcated in four
points
$x=\alpha_i,\, \beta_i,\,i=1,2$. 
We consider the local geometry of
$p:(\De_{1,1},\De_{1,1}\cap \widetilde C_\tau)\to \De$.
Let $\check \phi_\tau(x)$ be the polynomial of degree $4$ which is the
Taylor expansion of $\Phi_\tau(x)$ modulo $x^5$.
First we observe that this branched covering 
$p:(\De_{1,1},\De_{1,1}\cap \widetilde C_\tau)\to \De$
is topologically equivalent
to the one  where we replace $C_\tau$
by the curve $\check C_\tau=\{y=\check \phi_\tau(x)\}$.
Next, the situation for $\check C_\tau$ and its two flex tangents
 inside $\De_{1,1}$ is equivalent to the following explicit one:
\[
 C'_\tau: y=\xi_\tau(x),\quad \xi_\tau(x):=x^4-6\tau^2 x^2
\]
For this, we consider simply a homotopy
$\Xi_t(x)=t\,\xi_\tau(x)+(1-t)\check \phi_\tau(x)$. Except a finite number of 
$t=t_1,\dots,t_q$,  this family of curves  defines equivalent
covering over $\De$.
In this model $\xi_\tau(x)$, we have 
$\alpha_1,\,\alpha_2=\pm \tau$ and $\beta_1,\,\beta_2=\mp 3\tau$.
We choose $\{x=1\}$ as the fixed generic fiber i.e., $\eta_0=1$.
On the fiber $x=1$, 
$P,\,Q_1,Q_2$ are  the intersections of the line $x=1$ and $C,\, L_1,\,L_2$
respectively and
 we choose the generators as in Figure 7.
The base point $B$ is chosen on the circle $|y|=1$. The other
$d-1$ intersection points of
$C_\tau\cap \{x=1\}$ are outside of the unit disk and the generators
$g_1,\dots, g_{d-1}$ are
omitted in the figure.
The loops are oriented counterclockwise.
Now we consider the loop in the base space
$\ell\circ \omega\circ \ell\inv$ where
$\ell$ is the line segment from $x=1$ to $x=\beta_2+\eps',\,
\eps'<<(1-\beta_2)$ 
and $\omega$ is the loop $|\eta-\beta_2|=\eps'$.
It is now easy to see that the monodromy relation along this loop is
nothing but $S_1: \, g_dh_2g_d\inv h_2\inv$ as is expected. 

\begin{figure}[htb] 
\setlength{\unitlength}{1bp}
\begin{picture}(600,150)(-100,0)
 \put(-50,160){\special{epsfile=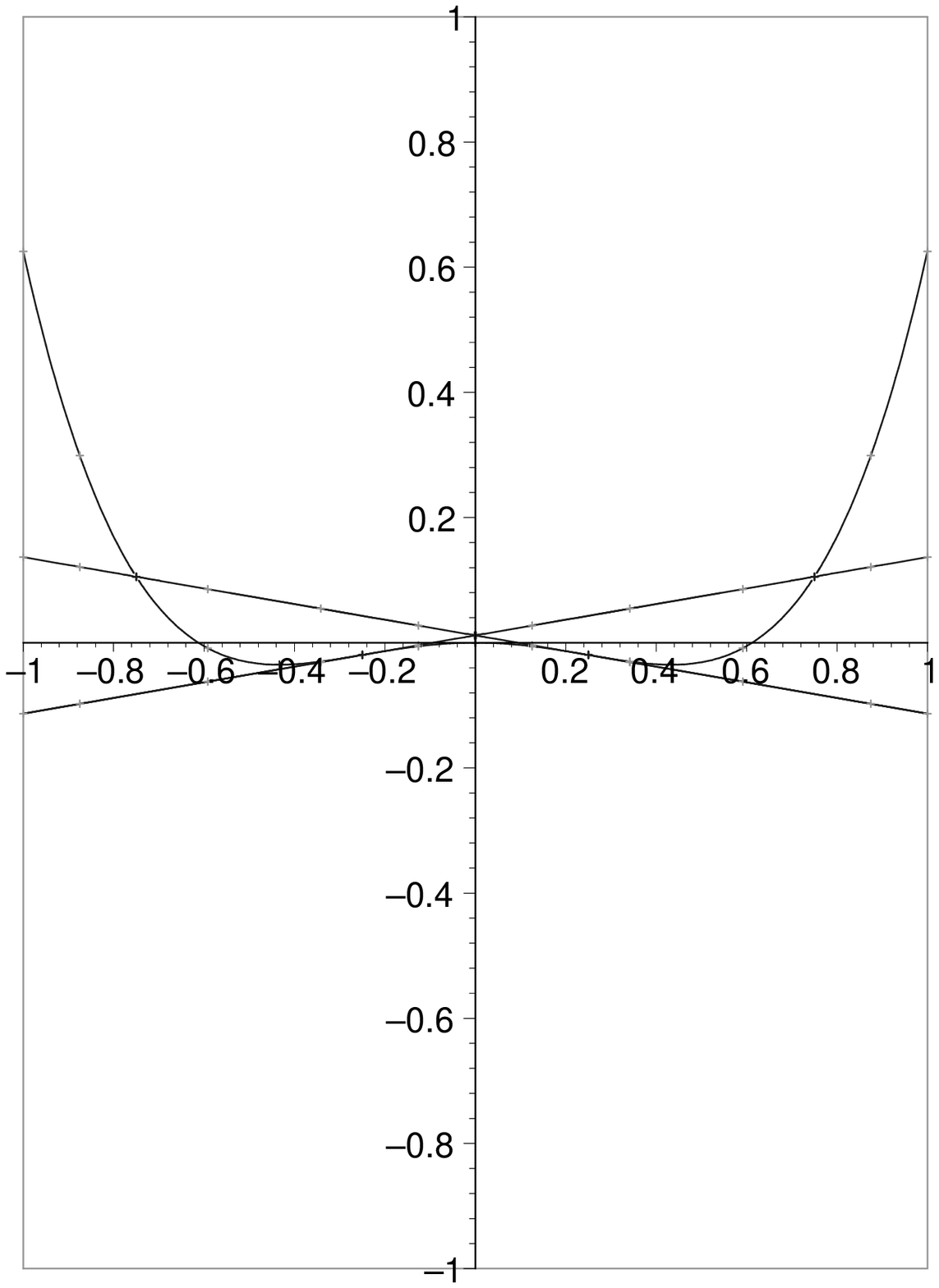 vscale=0.4
 hscale=0.4}}
 \put(140,130){\special{epsfile=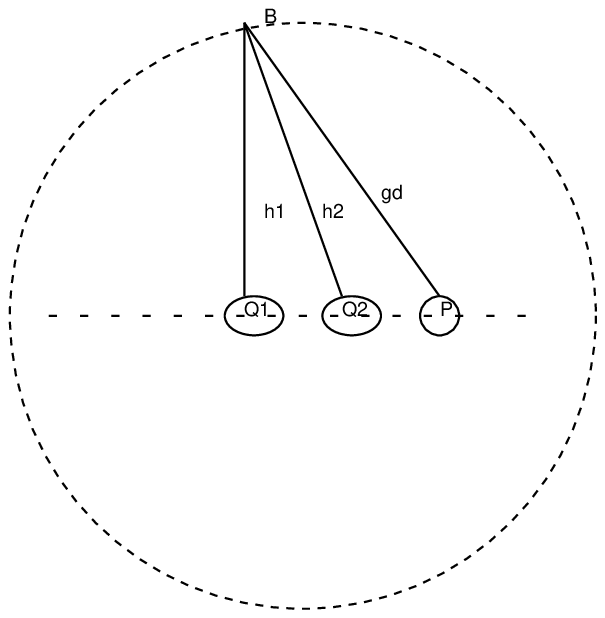 vscale=1.0
 hscale=1.0}}
\end{picture}
\vspace{1cm}
\caption{Choice of generators $g_d,h_1,h_2$}
\end{figure}
\noindent

\subsection{Graphs of various quartics}
We put the graphs  of various line-degenerated quartics of torus type.
\begin{figure}[H]
\setlength{\unitlength}{1bp}
\begin{picture}(600,150)(-100,0)
\put(-50,130){\special{epsfile=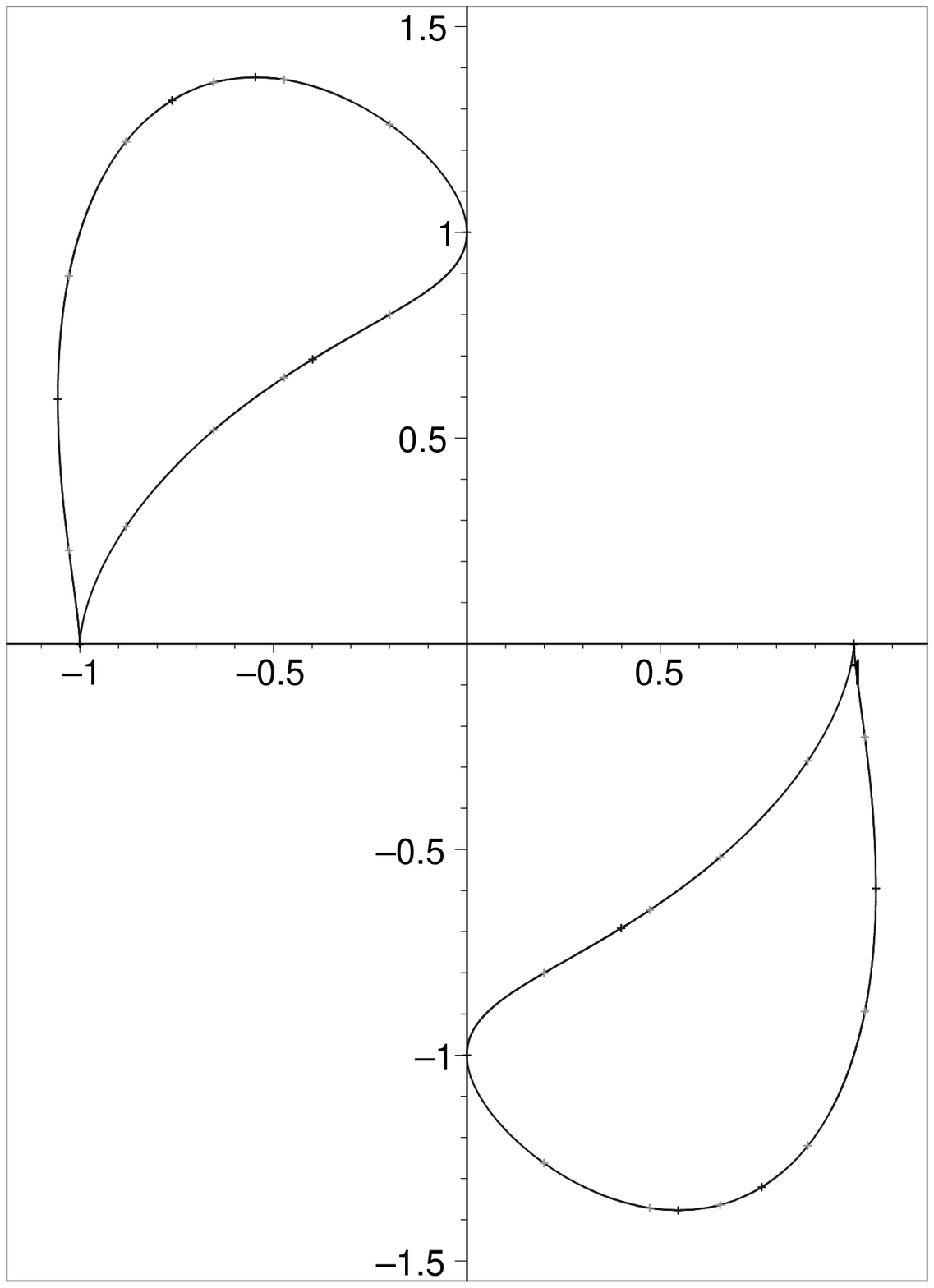 vscale=0.25
 hscale=0.3}}
\put(150,130){\special{epsfile=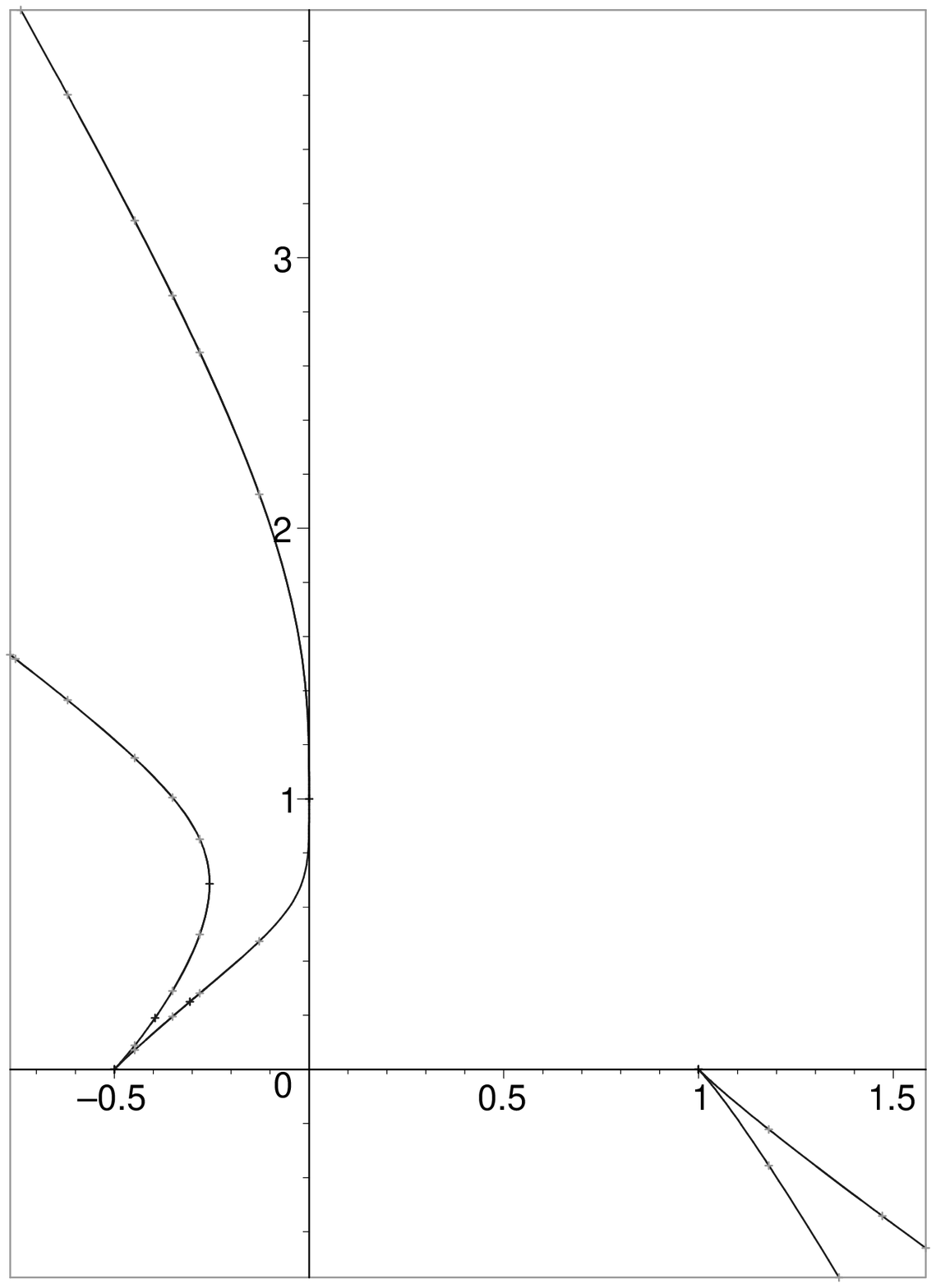 vscale=0.25
 hscale=0.3}}
\end{picture}
\vspace{0.4cm}
\caption{Quartic with $2A_2$, \,bi-tangent limit line (left),\,\, flex tangent limit line (right)}
\label{Quartic-2A2}
\end{figure}
\noindent

\begin{figure}[H]
\setlength{\unitlength}{1bp}
\begin{picture}(600,150)(-100,0)
\put(-50,130){\special{epsfile=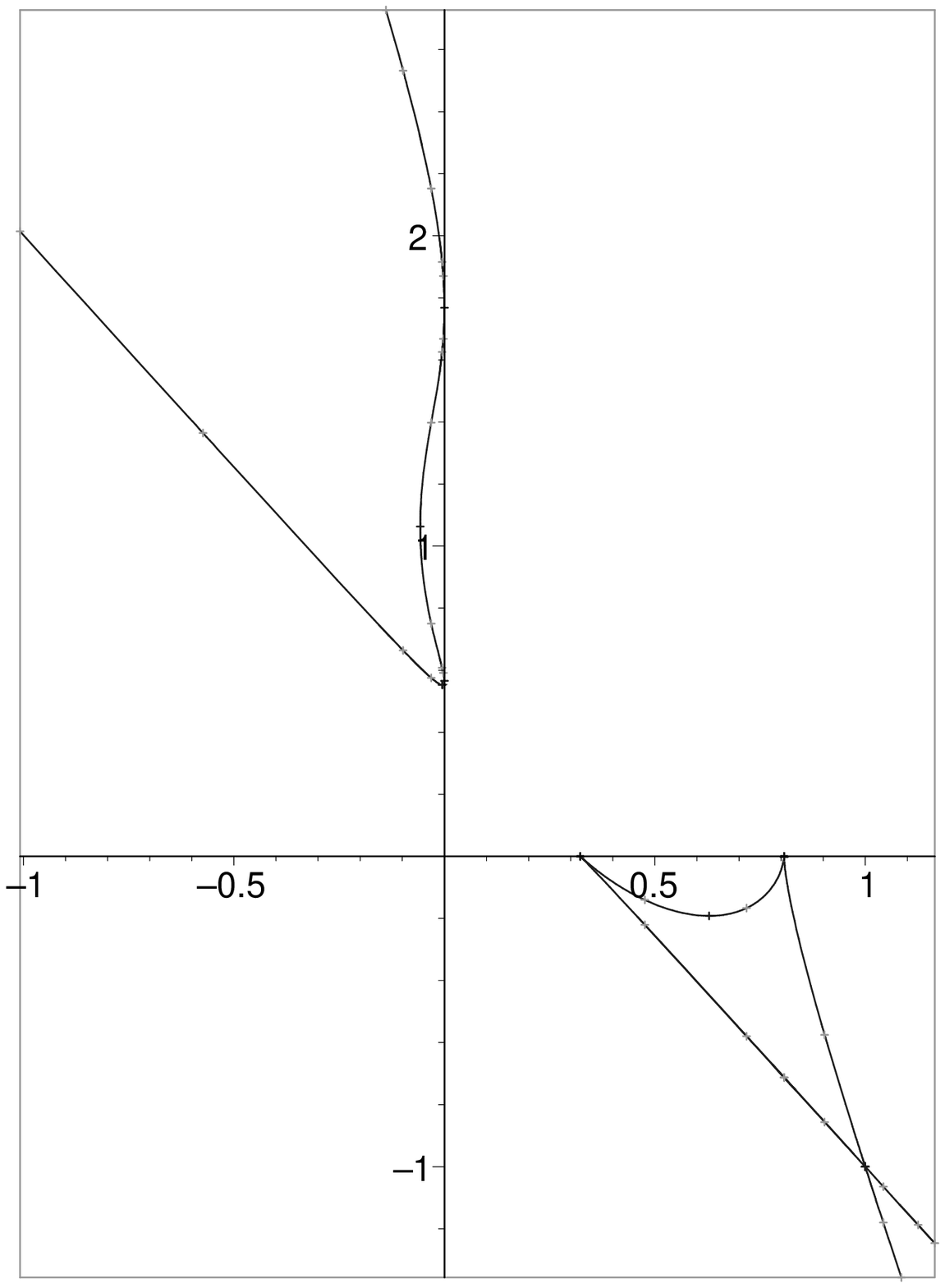 vscale=0.25
 hscale=0.3}}
\put(150,130){\special{epsfile=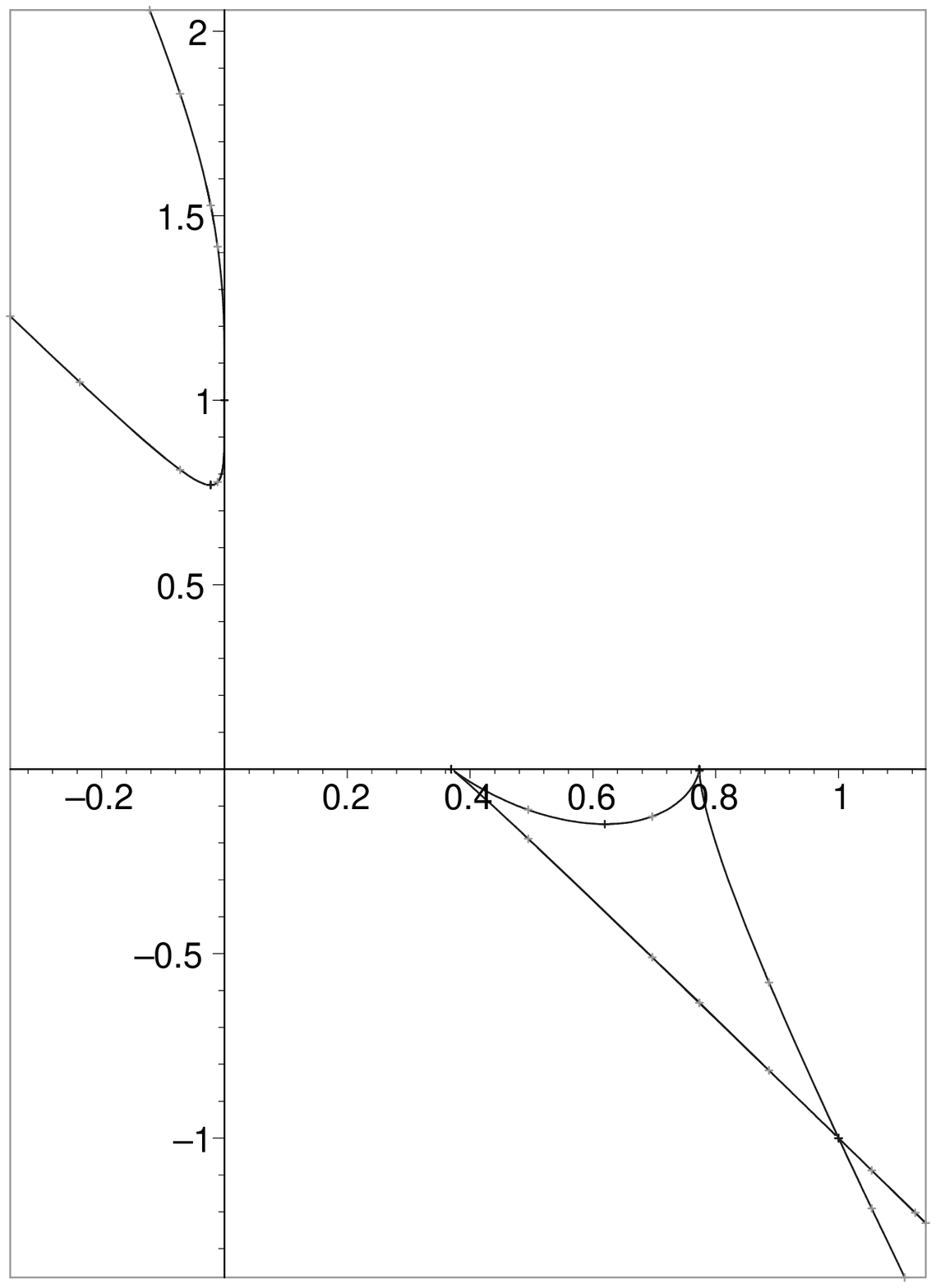 vscale=0.25
 hscale=0.3}}
\end{picture}
\vspace{.4cm}
\caption{Quartic with $2A_2+A_1$, \,bi-tangent limit line (left),\,\, flex tangent limit line (right)}
\label{Quartic-2A2A1}
\end{figure}
\noindent
\vspace{-1.cm}
\begin{figure}[H]
\setlength{\unitlength}{1bp}
\begin{picture}(600,150)(-100,0)
\put(-50,130){\special{epsfile=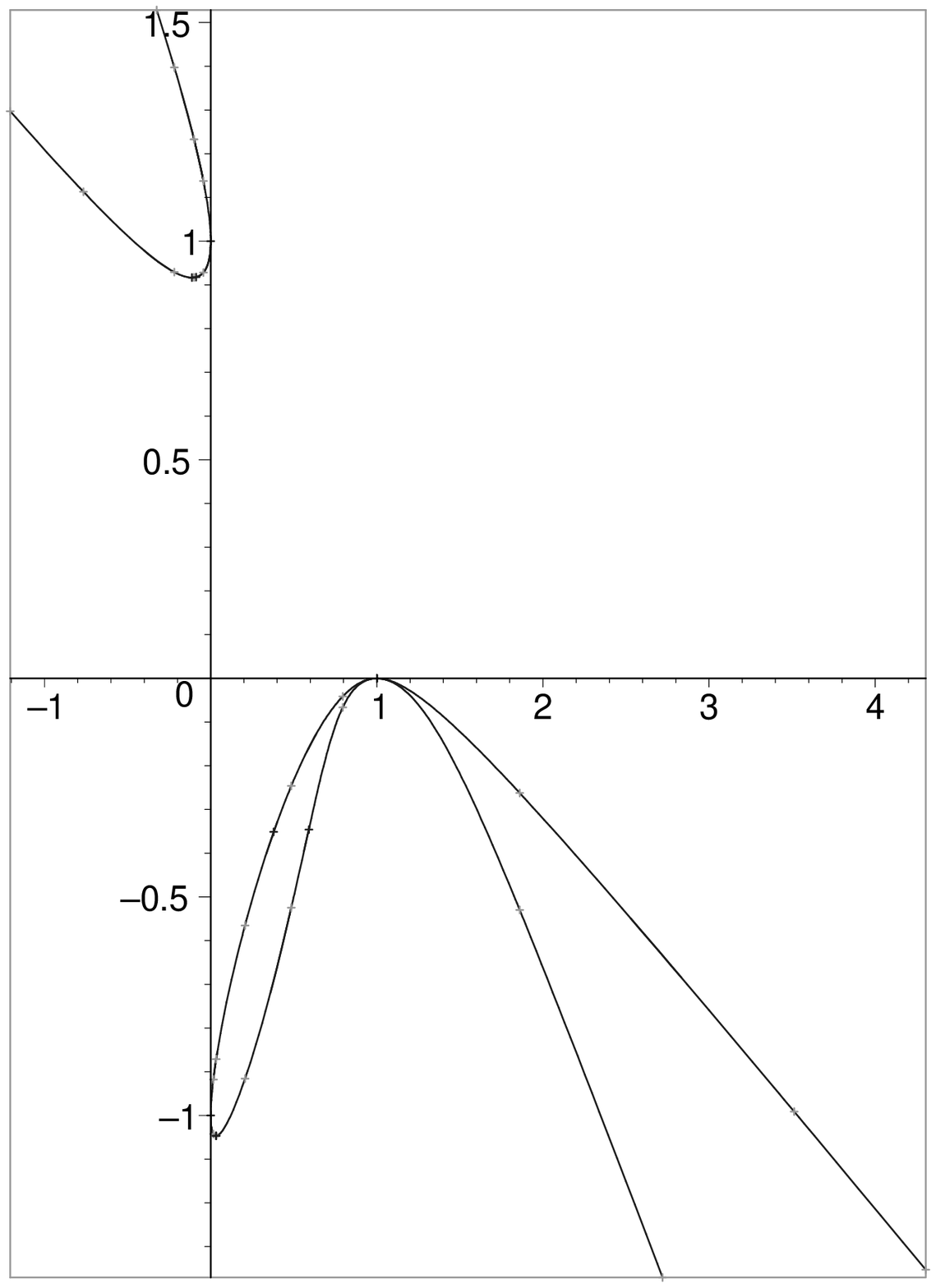 vscale=0.25
 hscale=0.3}}
\put(150,130){\special{epsfile=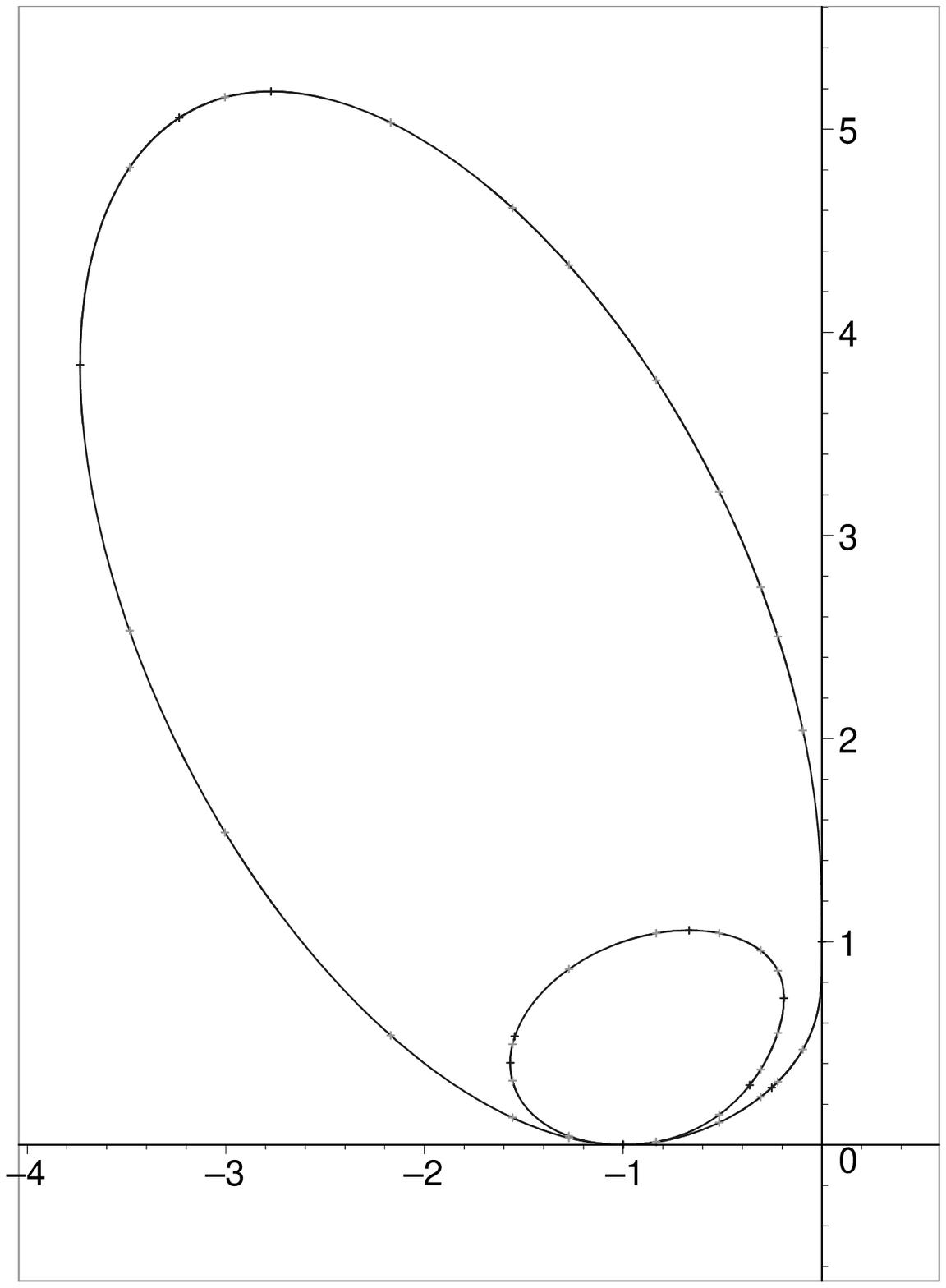 vscale=0.25
 hscale=0.3}}
\end{picture}
\vspace{.3cm}
\caption{Quartic with $A_5$, \,bi-tangent limit line (left),\,\, flex tangent
 limit line (right)}
\label{Quartic-A5}
\end{figure}
\vspace{-1cm}
\begin{figure}[H]
\setlength{\unitlength}{1bp}
\begin{picture}(600,150)(-100,0)
\put(-50,130){\special{epsfile=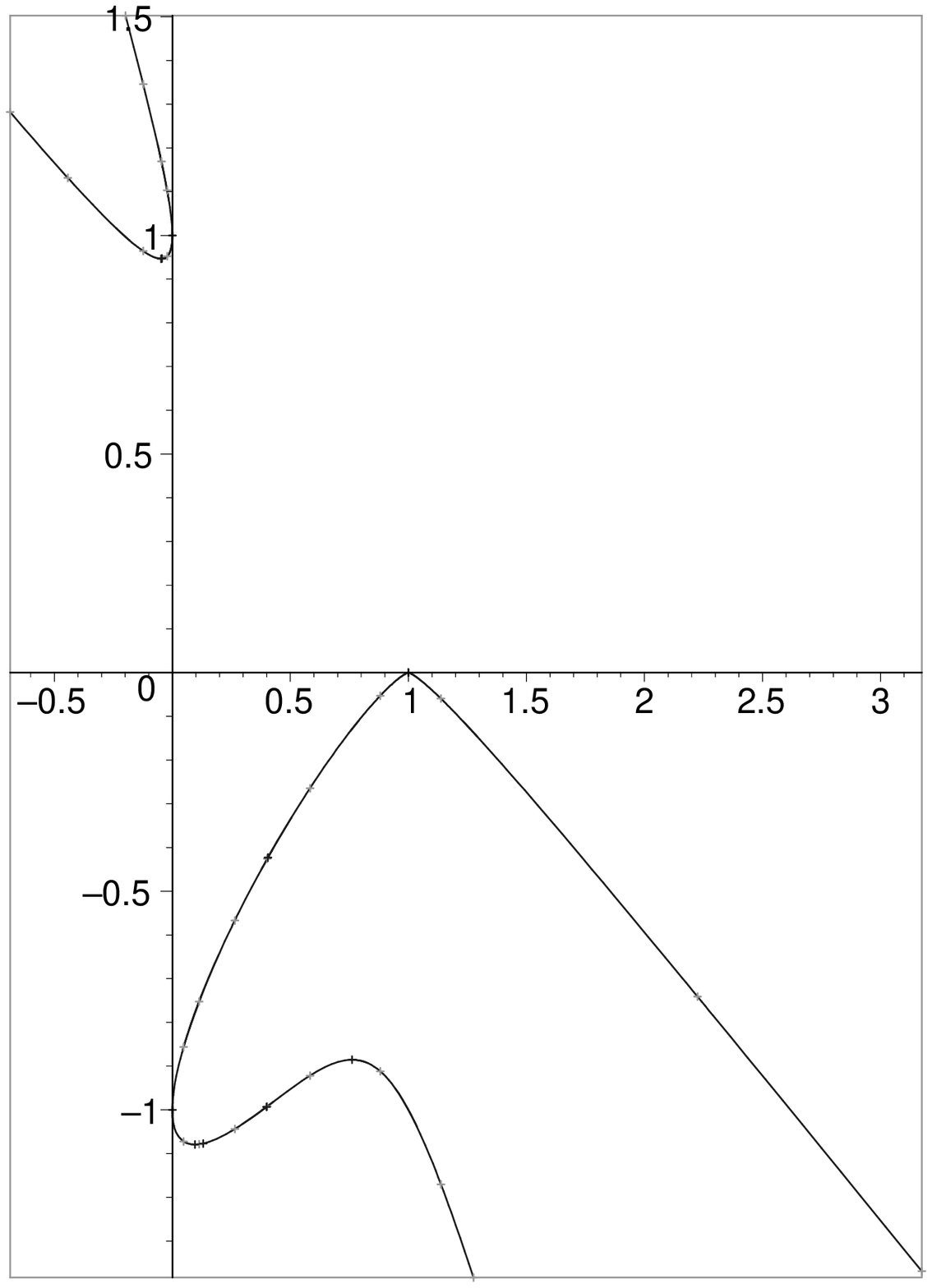 vscale=0.25
hscale=0.3}}
\put(150,130){\special{epsfile=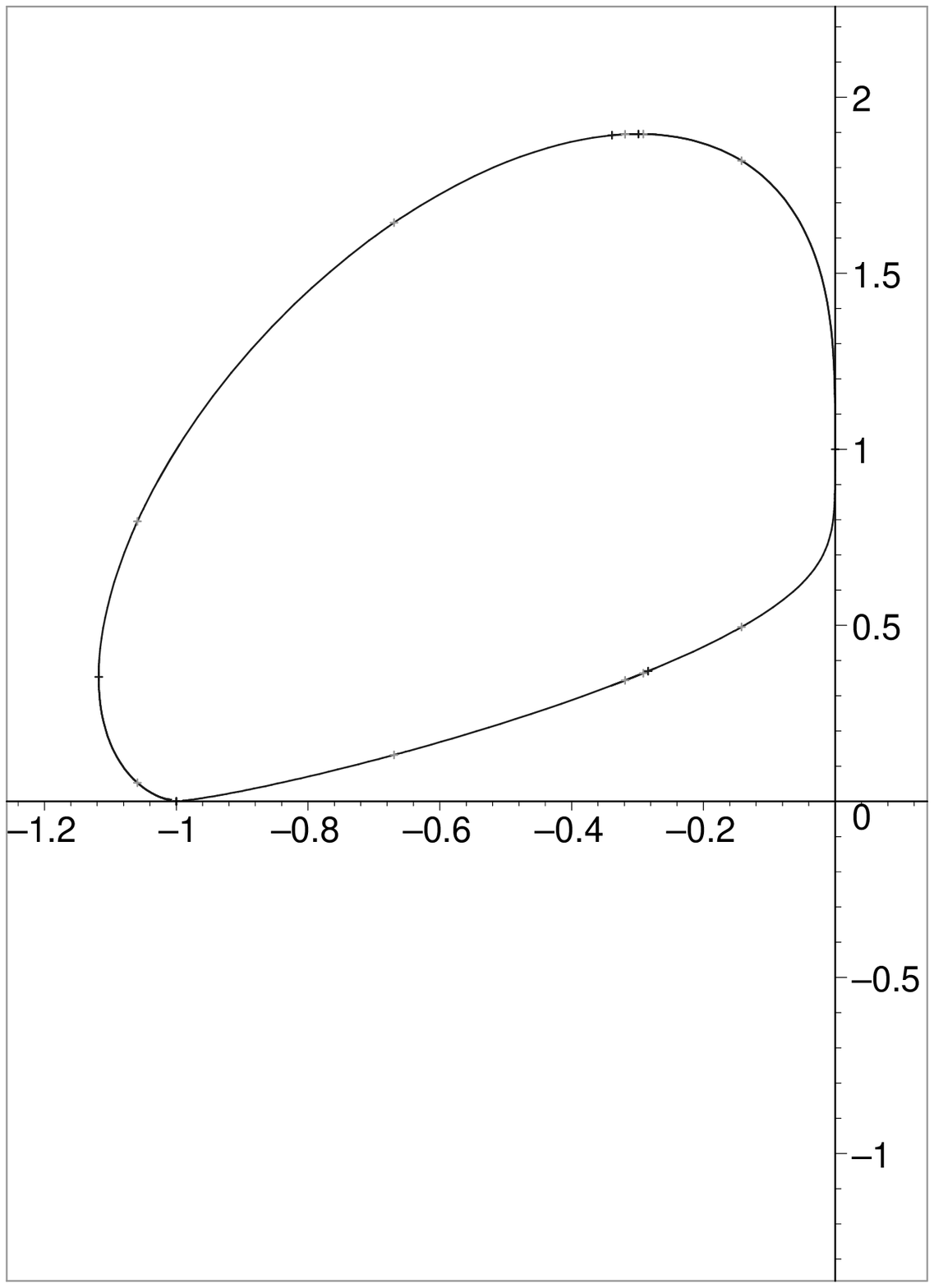 vscale=0.25
hscale=0.3}}
\end{picture}
\vspace{.3cm}
\caption{Quartic with $E_6$, \,bi-tangent limit line (left),\,\, flex tangent limit line (right)}
\label{Quartic-E6}
\end{figure}
\noindent
\begin{figure}[H]
\setlength{\unitlength}{1bp}
\begin{picture}(600,150)(-100,0)
\put(-50,130){\special{epsfile=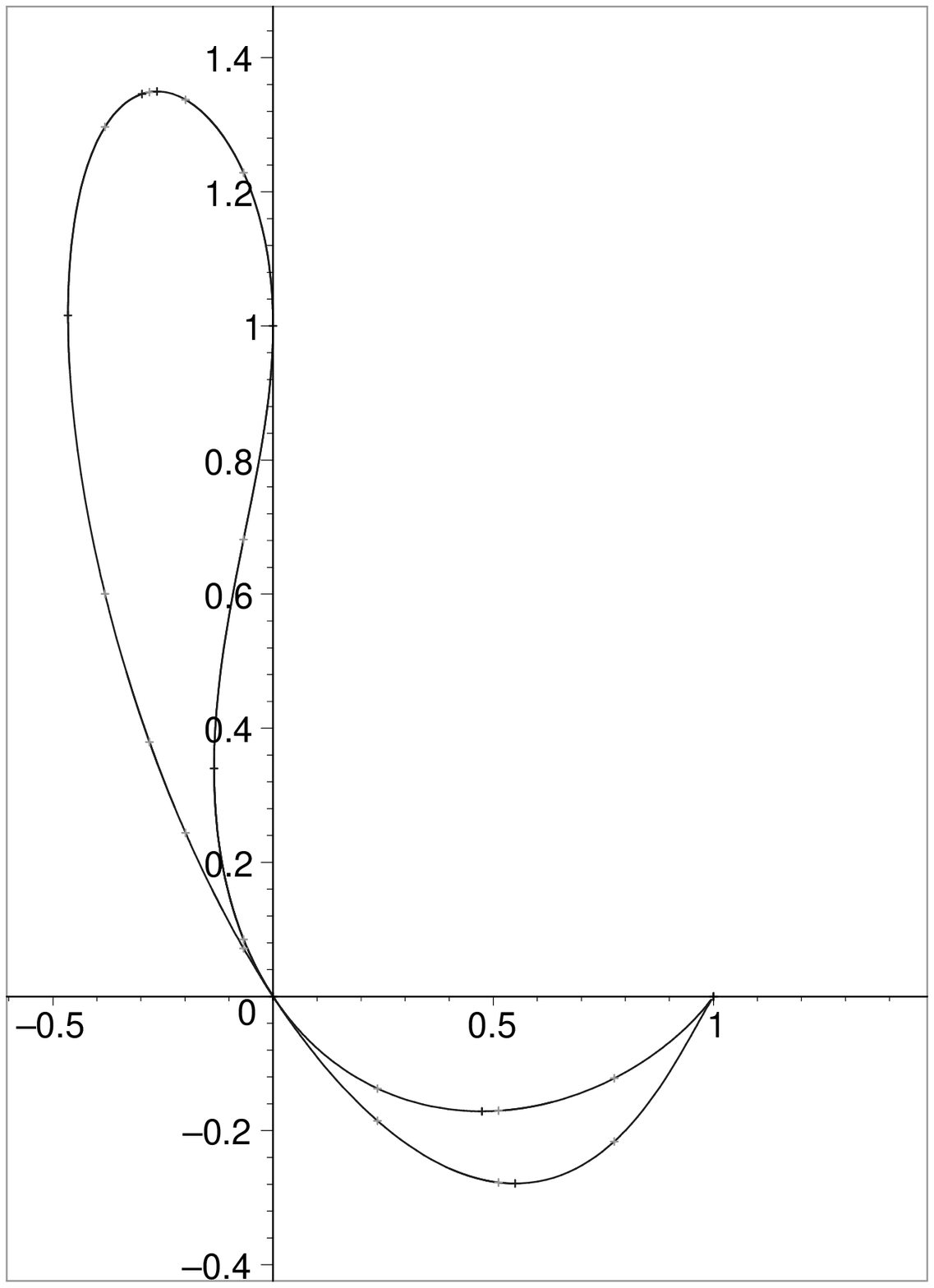 vscale=0.25
 hscale=0.3}}
\put(150,130){\special{epsfile=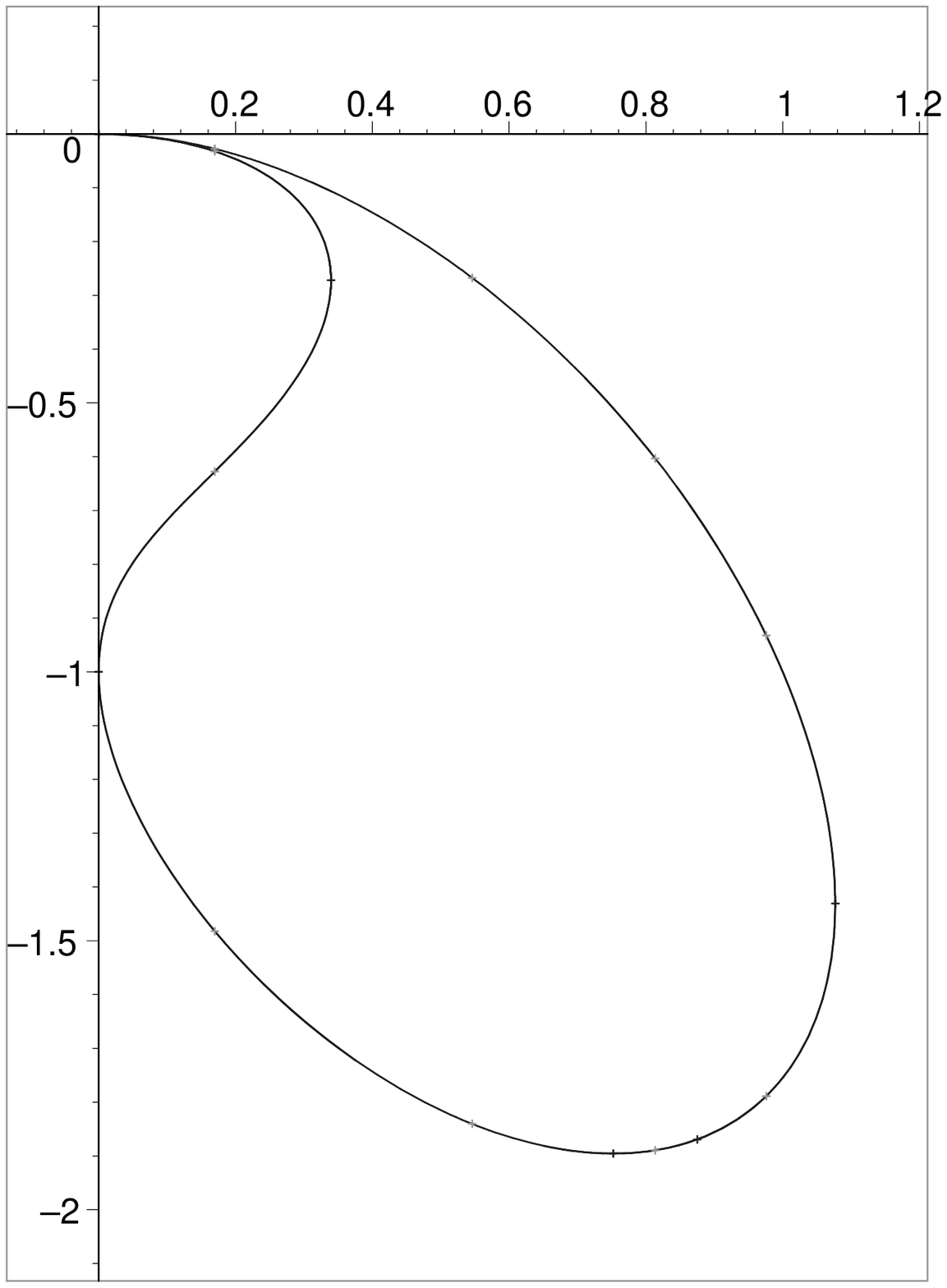 vscale=0.25
 hscale=0.3}}
\end{picture}
\vspace{.3cm}
\caption{Quartic with $A_3+A_2$\,\,left,\,\, with $A_6$\,\,right}
\label{Quartic-A3A2}
\end{figure}
\noindent

\def\cprime{$'$} \def\cprime{$'$} \def\cprime{$'$} \def\cprime{$'$}
  \def\cprime{$'$}

\end{document}